\documentclass[11pt]{article}
\usepackage{amssymb,amsmath,latexsym}
\usepackage{epsfig}

\newtheorem{thm}{Theorem}[section]
\newtheorem{cor}[thm]{Corollary}
\newtheorem{lemma}[thm]{Lemma}

\newtheorem{remark}[thm]{Remark}
\numberwithin{equation}{section}

\def\1{{\bf 1}}

\def\pf{{\medskip\noindent {\bf Proof. }}}
\def\qed{{\hfill $\Box$ \bigskip}}

\def\sA {{\cal A}}

  \def\sL {{\cal L}}

  \def\bR {{\mathbb R}}

\def\wt{\widetilde}

\def\E{{\mathbb E}}
\def\P{{\mathbb P}}

\def\bea{\begin{align*}}
\def\eea{\end{align*}}
\def\bee{\begin{equation}}
\def\eee{\end{equation}}

\def\eps{\varepsilon}

\makeatletter
\@addtoreset{equation}{section}

\makeatother

\begin{document}
\allowdisplaybreaks
\bibliographystyle{plain}

\title{\Large \bf
Heat Kernel Estimate for $\Delta+\Delta^{\alpha/2}$ in $C^{1,1}$ open sets}

\author{{\bf Zhen-Qing Chen}\thanks{Research partially supported
by NSF Grant DMS-0906743.}, \quad {\bf Panki Kim}\thanks{This work
was supported by Basic Science Research Program through the National
Research Foundation of Korea(NRF) grant funded by the Korea
government(MEST)(2009-0093131).} \quad and  \quad {\bf Renming Song}
}
\date{(January 29, 2010)}

\maketitle

\begin{abstract}
We consider a family of pseudo differential operators $\{\Delta+
a^\alpha \Delta^{\alpha/2}; \ a\in (0, 1]\}$ on $\bR^d$ for every
$d\geq 1$ that evolves continuously from $\Delta$ to $\Delta +
\Delta^{\alpha/2}$, where $\alpha \in (0, 2)$. It gives rise to a
family of L\'evy processes
$\{X^a, a\in (0, 1]\}$ in $\bR^d$, where $X^a$ is the sum of a
Brownian motion and an independent symmetric $\alpha$-stable process
with weight $a$. We establish sharp two-sided estimates for the heat
kernel of $\Delta + a^{\alpha} \Delta^{\alpha/2}$ with zero exterior
condition in a family of open subsets, including  bounded $C^{1, 1}$
(possibly disconnected) open sets. This heat kernel is also the
transition density of the sum of a Brownian motion and an
independent symmetric $\alpha$-stable process with weight $a$ in
such open sets.

Our result is the first sharp two-sided estimates for the transition
density of a Markov process with both diffusion and jump components
in open sets. Moreover, our result is uniform in $a$ in the sense
that the constants in the estimates are independent of $a\in (0, 1]$
so that it  recovers the Dirichlet heat kernel estimates for
Brownian motion by taking $a\to 0$. Integrating the heat kernel
estimates in time $t$, we recover the two-sided sharp uniform Green
function estimates of $X^a$ in bounded $C^{1,1}$ open sets in
$\bR^d$, which were recently established in \cite{CKSV2}   by using
a completely different approach.
\end{abstract}

\bigskip
\noindent {\bf AMS 2000 Mathematics Subject Classification}: Primary
60J35, 47G20, 60J75; Secondary
47D07

\bigskip\noindent
{\bf Keywords and Phrases}:
fractional Laplacian, Laplacian, symmetric
$\alpha$-stable process, Brownian motion,
heat kernel, transition density, L\'evy system
\bigskip

\section{Introduction}

Many physical and economic systems should be and in fact have been
successfully modeled by discontinuous Markov processes; see for
example, \cite{JW, KSZ, OS} and the references therein.
Discontinuous Markov processes are also very important from a
theoretical point of view, since they contain stable processes,
relativistic stable processes and jump diffusions as special cases.
Due to their importance both in theory and in applications,
discontinuous Markov processes have been receiving intensive study
in recent years.

In general, Markov processes may have both diffusion and jump
components. A Markov process having continuous sample paths is
called a diffusion. Diffusion processes in $\bR^d$ and second order
elliptic differential operators on $\bR^d$ are closely related in
the following sense. For a large class of second order elliptic
differential operators $\sL$ on $\bR^d$, there is a diffusion
process $X$ in $\bR^d$ associated with it so that $\sL$ is the
infinitesimal generator of $X$, and vice versa.  The connection
between $\sL$ and $X$ can also be seen as follows. The fundamental
solution $p(t, x, y)$ of $\partial _t u =\sL u$ (also called  the
heat kernel of $\sL$) is the transition density of $X$. Thus
obtaining sharp two-sided estimates for $p(t, x, y)$ is a
fundamental problem in both analysis and probability theory. In
fact, two-sided heat kernel estimates for diffusions in $\bR^d$ have
a long history and many beautiful results have been established. See
\cite{D1, D3} and the references therein. But, due to the
complication near the boundary, two-sided estimates on the
transition density of killed diffusions in a domain $D$
(equivalently, the Dirichlet heat kernel) have been established only
recently. See \cite{D2, D3, DS} for upper bound estimates and
\cite{Zq3} for the lower bound estimate of the Dirichlet heat
kernels in bounded $C^{1,1}$ domains.

The infinitesimal generator of a discontinuous Markov process in
$\bR^d$ is no longer a differential operator but rather a non-local
(or integro-differential) operator $\sL$. For instance, the
infinitesimal generator of a rotationally symmetric $\alpha$-stable
process in $\bR^d$ with $\alpha \in (0, 2)$ is a fractional
Laplacian operator $ \Delta^{\alpha /2}:=-  (-\Delta)^{\alpha/2}$.
Most of the recent studies concentrate on pure jump Markov
processes, like the rotationally symmetric $\alpha$-stable
processes, that do not have a diffusion component. For a summary of
some of these recent results from the probability literature, one
can see \cite{BBKRSV, C0} and the references therein. We refer the
readers to \cite{CSS, CaS, CV} for a sample of recent progresses in
the PDE literature.

Recently in \cite{CKS}, we obtained sharp two-sided estimates for
the heat kernel of the fractional Laplacian $\Delta^{\alpha /2}$ in
$D$ with zero exterior condition  (or equivalently, the transition
density function of the symmetric $\alpha$-stable process killed
upon exiting $D$) for any $C^{1, 1}$ open set $D\subset \bR^d$ with
$d\geq 1$. As far as we know, this was the first time
sharp two-sided estimates were established for Dirichlet heat
kernels of non-local operators. Since then, studies on this topic
have been growing rapidly. In \cite{CKS1, CKS2, CKS3}, the ideas of
\cite{CKS} were adapted to establish two-sided heat kernel estimates
of other pure jump Markov processes in open subsets of $\bR^d$. In
\cite{CT}, the large time behaviors of heat kernels for symmetric
$\alpha$-stable processes and censored stable processes in unbounded
open sets were studied. Very recently in \cite{BG, BGR}, the heat
kernel of the fractional Laplacian in non-smooth open set was
discussed. We refer the readers to \cite{C} for a survey on the
recent progresses in the heat kernel estimates of jump Markov
processes.

However until now, two-sided heat kernel estimates of Markov
processes with both diffusion and jump components in proper open
subsets of $\bR^d$  have not been studied. The fact that such a
process $X$ has both diffusion and jump components is the source of
many difficulties. The main difficulty stems from the fact that such
a process $X$ runs on two different scales: on the small scale the
diffusion part dominates, while on the large scale the jumps take
over. Another difficulty is encountered at the exit of $X$ from an
open set: for diffusions, the exit is through the boundary, while
for the pure jump processes, typically the exit happens by jumping
across the boundary. For a process $X$ that has both diffusion and
jump components, both cases will occur, which makes the process $X$
much more difficult to study.

In this paper, we consider L\'evy processes that are independent
sums of Brownian motions and (rotationally) symmetric stable processes
in $\bR^d$ with $d\geq 1$. We establish two-sided heat
kernel estimates for such L\'evy processes killed upon exiting a
$C^{1,1}$ open set. The processes studied in this paper serve as a
test case for more general processes with both diffusion and jump
components, just like Brownian motions do for more general
diffusions. We hope that our study will help to shed new light on
the understanding of the heat kernel behavior of more general Markov
processes. Although two-sided heat kernel estimates for Markov
processes with both diffusion and jump components in $\bR^d$ have
been studied recently in \cite{CK08, SV07}, as far as we know, this
is the first time that sharp two-sided estimates on the Dirichlet
heat kernels for Markov processes  with both  diffusion and jump
components in proper open subsets are established.

Let us now describe the main result of this paper and at the same
time fix the notations. Throughout this paper, we assume that $d\ge
1$ is an integer and $\alpha\in (0, 2)$. Let $X^0=(X^0_t,\, t\ge 0)$
be a Brownian motion in $\bR^d$ with generator $\Delta=\sum_{i=1}^d
\frac{\partial^2}{\partial x_i^2}$. Let $Y=(Y_t,\, t\ge 0)$ be a
(rotationally) symmetric $\alpha$-stable process in $\bR^d$, that
is, a L\'evy process such that
$$
\E_x \left[ e^{i\xi\cdot(Y_t-Y_0)} \right]\,=\,e^{-t|\xi|^{\alpha}}
\qquad \hbox{for every } x\in \bR^d \hbox{ and }  \xi\in \bR^d.
$$
The infinitesimal generator of a symmetric $\alpha$-stable process
$Y$ in $\bR^d$ is the fractional Laplacian $\Delta^{\alpha /2} $,
which is a prototype of non-local operators. The fractional
Laplacian can be written in the form
\begin{equation}\label{e:1.1}
\Delta^{\alpha /2} u(x)\, =\, \lim_{\eps \downarrow 0}\int_{\{y\in
\bR^d: \, |y-x|>\eps\}} (u(y)-u(x)) \frac{\sA (d, \alpha)
}{|x-y|^{d+\alpha}}\, dy,
\end{equation}
where $ {\cal A}(d, \alpha):= \alpha2^{\alpha-1}\pi^{-d/2}
\Gamma(\frac{d+\alpha}2) \Gamma(1-\frac{\alpha}2)^{-1}. $ Here
$\Gamma$ is the Gamma function defined by $\Gamma(\lambda):=
\int^{\infty}_0 t^{\lambda-1} e^{-t}dt$ for every $\lambda > 0$.
Assume that $X^0$ and $Y$ are independent. For any $a>0$, we define
$X^a$ by $X_t^a:=X^0_t+ a Y_t$. We will call the process $X^a$ the
independent sum of the Brownian motion $X^0$ and the symmetric
$\alpha$-stable process  $Y$ with weight $a>0$. The infinitesimal
generator of $X^a$ is $\Delta+ a^\alpha \Delta^{\alpha/2}$ and
$$
\E_x \left[ e^{i\xi\cdot(X^a_t-X^a_0)} \right]\,=\,e^{-t(|\xi|^{2}+a^{\alpha}|\xi|^{\alpha})}
\qquad \hbox{for every } x\in \bR^d \hbox{ and }  \xi\in \bR^d.
$$
Since
$$
a^{\alpha}|\xi|^{\alpha}= \int_{\bR^d}(1-\cos(\xi\cdot y))\,
\frac{ a^{\alpha} {\cal A}(d, \alpha)}{|y|^{d+\alpha}}dy,
$$
the density of the L\'evy measure of $X^a$ with respect
to the Lebesgue measure on $\bR^d$ is
$$
J^a(x, y)=j^a(|x-y|)
:=a^{\alpha}{\cal A}(d,
\alpha)|x-y|^{-(d+\alpha)}.
$$
The  function  $J^a (x, y)$  determines a L\'evy system for $X^a$,
which describes the jumps of the process $X^a$: for any non-negative
measurable function $f$ on $\bR_+ \times \bR^d\times \bR^d$ with
$f(s, y, y)=0$ for all $y\in \bR^d$, any stopping time $T$ (with
respect to the filtration of $X^a$) and any $x\in \bR^d$,
\begin{equation}\label{e:levy}
\E_x \left[\sum_{s\le T} f(s,X^a_{s-}, X^a_s)
\right]= \E_x \left[
\int_0^T \left( \int_{\bR^d} f(s,X^a_s, y) J^a(X^a_s,y) dy \right)
ds \right]
\end{equation}
(see, for example, \cite[Proof of Lemma 4.7]{CK} and \cite[Appendix
A]{CK2}). Let $p^a(t, x, y)$ be the transition density of the
process $X^a$ with respect to the Lebesgue measure on $\bR^d$, which
is known to exist and is jointly continuous on $(0, \infty)\times
\bR^d\times \bR^d$. For any $\lambda>0$, the process $(\lambda
X^{a}_{\lambda^{-2} t}, t\geq 0)$ has the same distribution as
$(X^{a \lambda^{(\alpha-2)/\alpha}}_t, t\ge 0)$ (see the second
paragraph of Section \ref{s:lb}), so we have
\begin{equation}\label{scalingrd}
p^{a\lambda^{(\alpha-2)/\alpha}} ( t,  x, y)= \lambda^{-d} p^{a}
(\lambda^{-2}t, \lambda^{-1} x, \lambda^{-1} y) \qquad \hbox{for }
t>0 \hbox{ and } x, y \in \bR^d.
\end{equation}

The following sharp two-sided estimates on $p^a(t, x, y)$ follows
from \eqref{scalingrd} and the main results in
\cite{CK08, SV07}
that give the sharp estimates on $p^1(t, x, y)$.

\begin{thm}\label{T:1.1}
There are constants $C_i\geq 1$, $i=1, 2$, such that, for all $a\in
[0, \infty)$ and $(t, x, y)\in (0, \infty]\times \bR^d\times \bR^d$
\begin{align*}
&C_1^{-1}  \,\left(t^{-d/2} \wedge (a^\alpha t)^{-d/\alpha} \right)
\wedge \left(t^{-d/2} e^{-C_2|x-y|^2/t}+ (a^{\alpha}
t)^{-d/\alpha}\, \wedge
\frac{a^{\alpha}t}{|x-y|^{d+\alpha}}  \right) \\
&\le \,p^a(t, x, y)\,\le\, C_1\, \,\left(t^{-d/2} \wedge
(a^{\alpha}t)^{-d/\alpha}\right) \wedge \left(t^{-d/2}
e^{-|x-y|^2/C_2t}+ (a^{\alpha} t)^{-d/\alpha}\, \wedge
\frac{a^{\alpha}t}{|x-y|^{d+\alpha}}   \right) .
\end{align*}
\end{thm}

Here for $a, b\in \bR$, $a\wedge b:=\min \{a, b\}$ and $a\vee
b:=\max\{a, b\}$. In particular, we have

\begin{cor}\label{C:1.2}
For any $M>0$ and $T>0$, there is a constant $C_3\geq 1$ depending
only on $M$ and $T$ such that, for all $a\in (0, M]$ and $(t, x,
y)\in (0, T]\times \bR^d\times \bR^d$
\begin{align*}
&C_3^{-1}  \, \left(t^{-d/2} e^{-C_2|x-y|^2/t}+ t^{-d/2}\, \wedge
\frac{a^{\alpha}t}{|x-y|^{d+\alpha}}   \right) \\
&\le \,p^a(t, x, y)\,\le\,C_3\, \left(  t^{-d/2}
e^{-|x-y|^2/C_2t}+ t^{-d/2}\, \wedge
\frac{a^{\alpha}t}{|x-y|^{d+\alpha}}  \right),
\end{align*}
where $C_2\geq 1$ is the constant in Theorem \ref{T:1.1}.
\end{cor}

Recall that an open set $D$ in $\bR^d$ (when $d\ge 2$) is said to be
$C^{1,1}$ if there exist a localization radius $ R>0 $ and a
constant $\Lambda>0$ such that for every $z\in\partial D$, there is
a $C^{1,1}$-function $\phi=\phi_z: \bR^{d-1}\to \bR$ satisfying
$\phi(0)=0$, $ \nabla\phi (0)=(0, \dots, 0)$, $\| \nabla \phi
\|_\infty \leq \Lambda$, $| \nabla \phi (x)-\nabla \phi (z)| \leq
\Lambda |x-z|$, and an orthonormal coordinate system $CS_z$:
$y=(y_1, \cdots, y_{d-1}, y_d):=(\wt y, \, y_d)$ with origin at $z$
such that $B(z, R )\cap D= \{y=({\tilde y}, y_d) \in B(0, R) \mbox{
in } CS_z: y_d > \phi (\wt y) \}$. The pair $( R, \Lambda)$ will be
called the $C^{1,1}$ characteristics of the open set $D$. By a
$C^{1, 1}$ open set in $\bR$ we mean an open set which can be
written as the union of disjoint intervals so that the minimum of
the lengths of all these intervals is positive and the minimum of
the distances between these intervals is positive. Note that a
$C^{1,1}$ open set can be unbounded and disconnected, and that a
bounded $C^{1,1}$ open set have only finitely many connected
components.

For an open set $D\subset \bR^d$ and $x\in D$, we will use
$\delta_D(x)$ to denote the Euclidean distance between $x$ and
$D^c$. For an open set $D\subset \bR^d$ and $(r_0, \lambda_0)\in
(0, \infty)\times [1, \infty)$, we say {\it the path distance in
each connected component of  $D$ is comparable to the Euclidean
distance with characteristics $(r_0, \lambda_0)$} if the following
holds for any $r\in (0, r_0)$: for every $x, y$ in the same
component of $D$ with $\delta_D(x) \wedge \delta_D(y)\geq r$, there
is a rectifiable curve $l$ in $D$ connecting $x$ to $y$ so that the
length of $l$ is no larger than $\lambda_0|x-y|$. Clearly, such a
property holds for all bounded $C^{1,1}$ open sets, $C^{1,1}$ open
sets with compact complements and domains above graphs of $C^{1,1}$
functions.

For any open subset $D\subset \bR^d$, we use $\tau^a_D$ to denote
the first time the process $X^a$ exits $D$. We define the process
$X^{a,D}$ by $X^{a,D}_t=X^a_t$ for $t<\tau^a_D$ and
$X^{a,D}_t=\partial$ for $t\ge \tau^a_D$, where $\partial$ is a
cemetery point. $X^{a,D}$ is called the subprocess of $X^a$ killed
upon exiting $D$. The infinitesimal generator of $X^{a,D}$ is
$(\Delta+ a^\alpha \Delta^{\alpha/2})|_D$. It follows from
\cite{CK08} that $X^{a,D}$ has a continuous transition density
$p^a_D(t, x, y)$ with respect to the Lebesgue measure.

The goal of this paper is to get the following sharp two-sided
estimates on $p^a_D(t, x, y)$ for any $C^{1,1}$ open set $D$ in
which the path distance in  each connected component of $D$ is
comparable to the Euclidean distance.

Let
\begin{equation}\label{eq:qd}
h^a_C(t, x, y):
=\begin{cases} &\left(1\wedge \frac{\delta_D(x)}{\sqrt{t}}
\right)\left(1\wedge \frac{\delta_D(y)}{\sqrt{t}}  \right) \left(
t^{-d/2} e^{-C|x-y|^2/t}+
\frac{a^\alpha t}{|x-y|^{d+\alpha}} \wedge t^{-d/2})\right)\\
&\qquad \qquad\hbox{when } x, y \hbox{ are in the same component of } D, \\
&\left(1\wedge \frac{\delta_D(x)}{\sqrt{t}}
\right)\left(1\wedge \frac{\delta_D(y)}{\sqrt{t}}  \right) \left (
 \frac{a^\alpha t}{|x-y|^{d+\alpha}}
\wedge t^{-d/2}\right)\\
&\qquad \qquad\hbox{when } x, y \hbox{ are
in different components of } D.
\end{cases}
\end{equation}

One can easily show that, when $D$ is bounded, the operator
$-(\Delta + a^\alpha \Delta^{\alpha/2})|_D$ has discrete spectrum
(see, for instance, the first paragraph of the proof of Theorem
\ref{t:main} (ii) and (iii) in Section 4).
In this case, we use $\lambda^{a, D}_1>0$ to denote the smallest
eigenvalue of $-(\Delta + a^\alpha \Delta^{\alpha/2})|_D$. Denote by
$D(x)$ the connected component of $D$ that contains $x$ and let $
\lambda^{a, D(x)}_1>0$ be the smallest eigenvalue of $-(\Delta +
a^\alpha \Delta^{\alpha/2})|_{D(x)}$.

\begin{thm}\label{t:main}
Let $d\geq 1$. Suppose that $D$ is a $C^{1,1}$
open set
in $\bR^d$
with characteristic $( R, \Lambda)$
such that the path distance
in  each connected component of $D$
is comparable to the
Euclidean distance
with characteristics $(r_0, \lambda_0)$.
\begin{description}
\item{\rm (i)}
For every $M>0$ and $T>0$, there are  positive constants $C_i=C_i(
R, \Lambda, r_0, \lambda_0, M, \alpha, T)\geq 1$, $i=4, 5$, such
that for all
$a \in (0, M]$ and $(t, x, y)\in (0, T]\times D\times
D$,
\begin{align*}
C_4^{-1}\, h^a_{C_5}(t, x, y) \le  p^a_D(t, x, y)
\leq C_4 h^a_{1/C_5}(t, x, y).
\end{align*}

\item{\rm (ii)}
Suppose in addition that $D$ is
bounded and connected.
For every $M>0$ and $T>0$,
there is a constant $C_6=C_6(D, M, \alpha, T)\geq 1$ so that for
all $a \in (0, M]$ and $(t, x, y)\in [T, \infty)\times D\times D$,
$$
C_6^{-1}\, e^{- t \, \lambda^{a, D}_1 }\, \delta_D (x)\, \delta_D (y)
\,\leq\,  p^a_D(t, x, y) \,\leq\, C_6\,
e^{- t\, \lambda^{a, D}_1}\, \delta_D (x) \,\delta_D (y).
$$

\item{\rm (iii)} Suppose   that $D$ is
bounded but disconnected. Then for every $M>0$ and $T>0$, there are
constants $C_i=C_i(D, M, \alpha, T)\geq 1$. $i=7,8$, such that for
all $a \in (0, M]$, $t\in [T, \infty)$, the following hold.

\noindent{\rm (a)}
If $x, y$ are in the same component $D(x)$ of $D$,
$$
C_7^{-1}\,
e^{- t \,  \lambda^{a, D(x)}_1 }
\, \delta_D (x)\, \delta_D (y)
\,\leq\,  p^a_D(t, x, y) \,\leq\,
C_7\left(e^{- t\, \lambda^{a, D(x)}_1 }+ \big(1\wedge (a^\alpha \, t)\big)
e^{- t\, \lambda^{a, D}_1 } \right) \delta_D (x) \,\delta_D (y).
$$
\noindent{\rm (b)}
If $x, y$ are in different components of $D$,
$$
 C_8^{-1}\, a^\alpha \, t \,
e^{- t \,  (\lambda^{a, D(x)}_1 \vee \lambda^{a, D(y)}_1) }\, \delta_D (x)\, \delta_D (y)
\,\leq\,  p^a_D(t, x, y) \,\leq\, C_8\,
(1 \wedge (a^\alpha  t) )\,
e^{- t\,\lambda^{a, D}_1 } \,
\delta_D (x) \,\delta_D (y).
$$
\end{description}
\end{thm}

\begin{remark}\label{R:1.4} \rm \begin{description}
\item{(i)}
Unlike the Brownian motion case, even though $D$ may be
disconnected, the process $X^{a, D}$ is always irreducible when
$a>0$ because $X^{a, D}$ can jump from one component of $D$ to
another. When $a>0$ is smaller, the connection between different
components of $D$ by $X^a$ becomes weaker. The estimates given in
Theorem \ref{t:main} present a precise quantitative description of
such a  phenomenon. Letting $a\to 0$, Theorem \ref{t:main} recovers
the Dirichlet heat kernel estimates for Brownian motion in
$D$ (even when $D$ is disconnected);
see
\cite{Ch, Zq3} and the reference therein for the latter. In
particular, for $x$ and $y$   in different components of $D$, we
have $\lim_{a\to 0+}p^a_D(t, x, y)=0$ for all $x, y>0$, which is the
case for Brownian motion.

\item{(ii)}
In fact, the estimates in Theorem \ref{t:main}(i) will be
established under a weaker assumption on $D$: the lower bounded
estimate is proved under the uniform interior ball condition and the
condition that the path distance in  each connected component of $D$
is comparable to the Euclidean distance (see Theorem \ref{t:low}),
while the upper bound estimate is proved under a weaker version of
the uniform exterior ball condition (see Theorem \ref{t:ub1}). Here
an open set $D\subset \bR^d$ is said to satisfy the {\it uniform
interior ball condition} with radius $R_1>0$ if for every $x\in D$
with $\delta_D(x)< R_1$, there is $z_x\in \partial D$ so that
$|x-z_x|=\delta_D(x)$ and $B(x_0, R_1)\subset D$ for
$x_0:=z_x+R_1(x-z_x)/|x-z_x|$. We say $D$ satisfies {\it a weaker
version of the  uniform exterior ball condition} with radius $R_1>0$
if for every $z\in \partial D$, there is a ball $B^z$ of radius
$R_1$ such that $B^z\subset (\overline D)^c$ and $\partial B^z  \cap
\partial D=\{z\}$.
\end{description}
\end{remark}

Integrating the  heat kernel estimates in Theorem \ref{t:main} over
time $t$ yields the following two-sided sharp estimates of the Green
function of $X^a$ in bounded $C^{1,1}$ open sets, which were first
obtained in \cite{CKSV2} by a different method. We will not give the
details in this paper on how these estimates can be obtained by
integrating the estimates in Theorem \ref{t:main}. Interested
readers are referred to the proof of \cite[Corollary 1.2]{CKS},
where the sharp estimates for the Green functions of symmetric
stable processes in bounded $C^{1,1}$ open sets are obtained from
the sharp heat kernel estimates for the heat kernels by integration
over time $t$.

\medskip

Define for $d\geq 3$ and $a>0$,
$$
g_D^a (x, y) := \begin{cases} \frac{1} {|x-y|^{d-2}} \left(1\wedge
\frac{  \delta_D(x) \delta_D(y)}{ |x-y|^{2}}\right)
\quad &\hbox{when } x, y \hbox{ are in the same component of } D, \\
\frac{a^\alpha} {|x-y|^{d-2}} \left(1\wedge \frac{  \delta_D(x)
\delta_D(y)}{ |x-y|^{2}}\right) \quad &\hbox{when } x, y \hbox{ are
in different components of } D;
\end{cases}
$$
for $d=2$ and $a>0$,
$$
g_D^a (x, y): = \begin{cases} \log\left(1+\frac{  \delta_D(x)
\delta_D(y)}{ |x-y|^{2}}\right)
\quad &\hbox{when } x, y \hbox{ are in the same component of } D, \\
a^\alpha \log\left(1+\frac{  \delta_D(x) \delta_D(y)}{
|x-y|^{2}}\right) \quad &\hbox{when } x, y \hbox{ are in different
components of } D;
\end{cases}
$$
and for $d=1$ and $a>0$,
$$
g_D^a (x, y): =
\begin{cases} \left(\delta_D(x)
\delta_D(y)\right)^{1/2}\wedge\frac{ \delta_D(x) \delta_D(y)}{|x-y|}
\quad &\hbox{when } x, y \hbox{ are in the same component of } D, \\
a^\alpha \big(\left(\delta_D(x) \delta_D(y)\right)^{1/2}\wedge\frac{
\delta_D(x) \delta_D(y)}{|x-y|}\big)\quad &\hbox{when } x, y \hbox{
are in different components of } D.
\end{cases}
$$

\begin{cor}\label{C:1.5}
Let $M>0$. Suppose that $D$ is a bounded $C^{1,1}$
open set in $\bR^d$. There exists
$C_9=C_9(D, M, \alpha)>1$ such that for
all $x, y \in D$ and all $a \in (0,M]$
$$
C_9^{-1} \, g_D^a(x, y) \leq G_D^a(x, y) \leq C_9\, g_D^a(x, y) .
$$
\end{cor}

This paper is a natural continuation of
\cite{CKS}, where sharp two-sided heat kernel estimates for symmetric
$\alpha$-stable processes in $C^{1,1}$ open sets are first derived,
as well as  \cite{CKSV1}, where the boundary Harnack principle for $X^a$ is
established. Some ideas of the approach in this paper can be traced
back to \cite{CKS} but a number of new ideas are needed to handle
the combined effects of Brownian motion
and discontinuous
stable process. A
comparison with  subordination of killed Brownian motion is used for
the lower bound
short time  heat kernel estimates for $X^{a, D}$.
We would like to point out that, unlike
\cite{CKS}, the boundary Harnack principle for $X^a$ is not used
directly in this paper. Instead we use one of the key lemmas
established in \cite{CKSV1} to obtain the upper bound of the heat
kernel (see Lemma \ref{L:2.0_1}). Theorem \ref{t:main}(i) will be
established through Theorem \ref{t:low} and Theorem \ref{t:ub1},
which give the lower bound and upper bound estimates, respectively.
In contrast to  that in \cite{CKS, CKS1, CKS2, CKS3}, the proof of
large time heat kernel estimates in Theorem \ref{t:main}(ii)-(iii)
does not use intrinsic ultracontractivity of
$X^{a,D}$.
The proof presented here is more direct, and uses only the continuity of
$\lambda^{a, D}_1$ and its corresponding first eigenfunction in
$a\in (0, M]$, which is established in \cite{CS8}. Lastly, we point
out  that the approach of \cite{BGR} relies critically on the fact
the symmetric stable processes do not have diffusion component and
so it is not directly applicable to the processes considered in this
paper.

\medskip

We will use capital letters $C_1, C_2, \dots$ to denote constants in
the statements of results, and their labeling will be fixed. The
lower case constants $c_1, c_2, \dots$ will denote generic constants
used in proofs, whose exact values are not important and can change
from one appearance to another. The labeling of the lower case
constants starts anew in each proof. The dependence of the constant
$c$ on the dimension $d$ will not be mentioned explicitly. We will
use ``$:=$" to denote a definition, which is read as ``is defined to
be". We will use $\partial$ to denote a cemetery point and for every
function $f$, we extend its definition to $\partial$ by setting
$f(\partial )=0$. We will use $dx$ to denote the Lebesgue measure in
$\bR^d$. The Lebesgue measure of a Borel set $A\subset \bR^d$ will
be denoted by $|A|$.

\section{Lower bound estimate}\label{s:lb}

In this section, we assume that $D$ is an open set in $\bR^d$
satisfying the {\it  uniform interior ball condition} with radius
$R_1>0$  and that the path distance in each connected component of
$D$ is comparable to the Euclidean distance with characteristics
$(r_0, \lambda_0)$. Observe that under the uniform interior ball
condition, the condition that the path distance for each connected
component in $D$ is comparable to the Euclidean distance is
equivalent to the following: there exist $r_2, \lambda >0$ such that
for all $r\in (0, r_2]$ and all $x,y$   in the same connected
component of $D$ with $\delta_D(x) \wedge \delta_D(y)\geq r$, there
is a rectifiable curve $l$ in $D$ connecting $x$ to $y$ so that the
length of $l$ is no larger than $ \lambda |x-y|$ and $\delta_D(z)
\geq r$ for every $z\in l$. The latter  is also equivalent to the
following, which is called the connected ball condition in
\cite{Ch}: For all $r\in (0, r_2]$ and $x,y$   in the same connected
component of $D$ with $\delta_D(x) \wedge \delta_D(y)>r$, there
exist $ m $ and $x_k$, $k= 1, 2, \ldots , m$ such that $x_0 = x, $
$x_m = y ,$ $ x_{k-1} \in B(x_k, \frac{r}{2}) \subset B( x_k , r)
\subset D$ and  $r \cdot m \leq \lambda_0 |x-y|$.

Observe for all $\lambda, a>0$ and $\xi, x \in \bR^d$,
\begin{align*}
\E_x\left[e^{i\xi\cdot(\lambda (X_{t/\lambda^{2}}^{a}-
X_0^{a}))}\right] = e^{-t|\xi|^2} \, \E_x\left[e^{i(a\lambda\xi)
\cdot(Y_{t/\lambda^{2}}- Y_0)}\right]
=e^{-t(|\xi|^2+(a\lambda^{(\alpha-2)/\alpha})^\alpha|\xi|^\alpha)}.
\end{align*}
It follows that
if $\{X^{a, D}_t, t\geq 0\}$ is the  subprocess in $D$ of the
independent sum of a Brownian motion and a symmetric $\alpha$-stable
process on $\bR^d$ with weight $a$, then $\big\{\lambda X^{a,
D}_{\lambda^{-2} t}, t\geq 0\big\}$ is the subprocess in $\lambda D$ of
the independent sum of a Brownian motion and a symmetric
$\alpha$-stable process on $\bR^d$  with  weight $a
\lambda^{(\alpha-2)/\alpha}$. So for any $\lambda>0$, we have
\begin{equation}\label{e:scaling}
p^{a\lambda^{(\alpha-2)/\alpha}}_{\lambda D} ( t,  x, y)=
\lambda^{-d} p^{a}_D (\lambda^{-2}t, \lambda^{-1} x, \lambda^{-1} y)
\qquad \hbox{for } t>0 \hbox{ and } x, y \in \lambda D.
\end{equation}
The above scaling property of $X^{a}$ will be used throughout this
paper. For $t >0$, we define
\begin{equation}\label{e:2.2}
a_t:= a t^{(2-\alpha)/(2\alpha)}.
\end{equation}
This
notation will be used in this paper when we scale an open $D$ by
$s^{-1/2}$ to $s^{-1/2}D$.

We first recall the definition of subordinate killed Brownian
motion: Assume that $U$ is an open subset in $\bR^d$ and $T_t$ is an
$\alpha/2$-stable subordinator independent of the killed Brownian
motion $X^{0, U}$. For each $a\ge0$, let $T^a$ be the subordinator
defined by $T^a_t:=t+a^2T_t$. Then the process $\{Z^{a, U}_t: t\ge
0\}$ defined by $Z^{a, U}_t=X^{0, U}_{T^a_t}$ is called a
subordinate killed Brownian motion in $U$. Let $q^{a}_U(t, x, y)$ be
the transition density of $Z^{a, U}$. Then it follows from
\cite[Proposition 3.1]{SV08} that
\begin{equation}\label{e:pq}
p^{a}_{U}(t, z, w)\ge q^{a}_{U}(t, z, w), \quad (t, z, w)\in (0,
\infty)\times U\times U.
\end{equation}
We will use this fact in the next result.

\begin{lemma}\label{lbbyskbm}
Suppose that $M$ and $T$ are positive constants.  Then there exist
positive constants
$C_i= C_i( R_1, r_0, \lambda_0,
\alpha, T, M)$, $i=10, 11$,
such that for all $a\in (0, M]$,  $t \in (0, T]$ and
$x, y$ in the same connected component of $D$,
$$
p^a_D(t, x, y)\ge C_{10} t^{-d/2}
\left(1\wedge \frac{\delta_D(x)}{\sqrt{t}}
\right)\left(1\wedge \frac{\delta_D(y)}{\sqrt{t}}
\right)e^{-C_{11}|x-y|^2/t} .
$$
\end{lemma}

\pf Suppose that $x$ and $y$ are in the same component, say $U$, of
$D$. Let $p_U(t, x, y)$ be the transition density of the killed
Brownian motion in $U$. It follows from \cite[Theorem 3.3]{Ch} (see
also \cite[Theorem 1.2]{Zq3}) that there exist positive constants
$c_1=c_1(R_1, r_0, \lambda_0, \alpha, T)$ and $c_2=c_2(R_1, r_0,
\lambda_0, \alpha)$
such that for any $(s, x, y)\in (0, 2T]
\times U\times
U$,
$$
p_U(s, x, y)\ge c_1\left(1\wedge \frac{\delta_U(x)}{\sqrt{s}}
 \right)\left(1\wedge \frac{\delta_U(y)}{\sqrt{s}}
\right)s^{-d/2}e^{-c_2|x-y|^2/s}.
$$
(Although not
explicitly mentioned in \cite{Ch}, a careful exam of the proofs in
\cite{Ch} reveals that the constants $c_1$ and $c_2$ in above lower
bound estimate can be chosen to depend only on $(R_1, r_0,
\lambda_0, \alpha, T)$ and  $(R_1, r_0, \lambda_0,
\alpha)$, respectively.)
Since $p_{t^{-1/2}U}(u, t^{-1/2}x, t^{-1/2}y)=t^{d/2}p_U(ut, x, y)$,
we have for $t \le T$
and $(u, x, y)\in (0, 2]\times U\times U$,
\begin{equation}\label{e:2.4}
p_{t^{-1/2}U}(u, t^{-1/2}x, t^{-1/2}y) \geq
c_1\left(1\wedge \frac{\delta_U(x)}{\sqrt{tu}}
\right)\left(1\wedge \frac{\delta_U(y)}{\sqrt{tu}}
\right)u^{-d/2}e^{-c_2|x-y|^2/(tu)}.
\end{equation}

Let $\mu^{a_t}(u, s)$ be the
density of ${a_t}^2T_u$,
where $a_t$ is defined in \eqref{e:2.2}.
Then it follows from the definition of
the subordinate killed Brownian motion
(for example, see \cite[page 149]{BBKRSV})
that
for every $1/3 \le b \le 1$ and $0 < t \le T$,
\begin{eqnarray*}
q^{a_t}_{t^{-1/2}U}(b, t^{-1/2}x, t^{-1/2}y)
&=&\int^{\infty}_bp_{t^{-1/2}U}(s, t^{-1/2}x, t^{-1/2}y)\P(b+a_t^2T_b \in ds)\\
&=&\int^{\infty}_bp_{t^{-1/2}U}(s, t^{-1/2}x, t^{-1/2}y)\mu^{a_t}(b, s-b)ds\\
&=&\int^{\infty}_0p_{t^{-1/2}U}(s+b, t^{-1/2}x, t^{-1/2}y)\mu^{a_t}(b, s)ds.
\end{eqnarray*}
Consequently, by  \eqref{e:pq} and \eqref{e:2.4},
for every $1/3 \le b \le 1$ and $0 < t \le T$,
\begin{eqnarray}
&& p^{a_t}_{t^{-1/2}D}(b, t^{-1/2}x, t^{-1/2}y) \nonumber \\
& \ge &
p^{a_t}_{t^{-1/2}U}(b, t^{-1/2}x, t^{-1/2}y) \nonumber\\
& \ge \, & q^{a_t}_{t^{-1/2}U}
(b, t^{-1/2}x, t^{-1/2}y) \nonumber\\
&\ge& \int^1_0p_{t^{-1/2}U}(s+b, t^{-1/2}x, t^{-1/2}y)
             \mu^{a_t}(b, s)ds\nonumber\\
& \ge  & \frac{c_1}{2}
\left(1\wedge \frac{\delta_U(x)}{\sqrt{t}}  \right)
\left(1\wedge \frac{\delta_U(y)}{\sqrt{t}}  \right) e^{-3c_2|x-y|^2/t}
\int^1_0\mu^{a_t}(b, s)ds\nonumber\\
& =&  \frac{c_1}{2}\left(1\wedge \frac{\delta_U(x)}{\sqrt{t}}  \right)
\left(1\wedge \frac{\delta_U(y)}{\sqrt{t}}  \right)
e^{-3c_2|x-y|^2/t} \,
\P({a_t}^2T_b\le 1)\nonumber\\
& \ge &
\frac{c_1}{2}\left(1\wedge \frac{\delta_U(x)}{\sqrt{t}}  \right)
\left(1\wedge \frac{\delta_U(y)}{\sqrt{t}}  \right)
e^{-3c_2|x-y|^2/t}\,
\P(T_{1/3}\le M^{-2}T^{-(2-\alpha)/\alpha})\nonumber\\
& \ge & c_3\left(1\wedge \frac{\delta_U(x)}{\sqrt{t}} \right)
\left(1\wedge \frac{\delta_U(y)}{\sqrt{t}}  \right)
e^{-3c_2|x-y|^2/t}\nonumber\\
&=& c_3\left(1\wedge \frac{\delta_D(x)}{\sqrt{t}} \right)
\left(1\wedge \frac{\delta_D(y)}{\sqrt{t}}  \right)
e^{-3c_2|x-y|^2/t}. \label{e:direct}
\end{eqnarray}

We now conclude from
\eqref{e:scaling}, \eqref{e:2.2}
and \eqref{e:direct} with $b=1$
that
$$
p^a_D(t, x, y)=t^{-d/2}p^{a_t}_{t^{-1/2}D}(1, t^{-1/2}x, t^{-1/2}y)
\ge
c_3t^{-d/2}\left(1\wedge \frac{\delta_D(x)}{\sqrt{t}}   \right)
\left(1\wedge \frac{\delta_D(y)}{\sqrt{t}}   \right)
e^{-3c_2|x-y|^2/t}.
$$
\qed

The inequality
\eqref{e:direct} above with $b=1/3$ will be used later.

\begin{lemma}\label{L:2.2}
For all $M, r, b >0$, there exists $C_{12}=C_{12}(M,r,b)>0$ such that
$$
\P_0(\tau^{a}_{B(0,r)}>b) \ge C_{12} >0 \qquad \hbox{for all } a \in
(0, M].
$$
\end{lemma}

\pf It follows from Lemma \ref{lbbyskbm} that
\begin{eqnarray*}
\inf_{a\in (0, M]}\P_0(\tau^{a}_{B(0,r)}>b)
&=&\inf_{a\in (0, M]}\int_{B(0,r)} p^a_{B(0,r)} (b,0,y) dy \\
&\ge& c \, b^{-d/2}\left(\frac{r}{\sqrt{b}} \wedge 1
\right)\int_{B(0,r)} \left(\frac{|y|}{\sqrt{b}} \wedge 1 \right)
e^{-c_1|y|^2/b} dy.\end{eqnarray*} \qed

\begin{lemma}\label{lower bound12}
Suppose that $M$ and $r$ are
positive constants.  Then there is a
constant $C_{13} = C_{13} (M, \alpha, d, r)\in (0, 1/3]$ such that
for all $a \in (0, M]$ and $u, v\in \bR^d$,
\begin{eqnarray*}
p^{a}_{B(u,r)\cup B(v,r)}(1/3, u, v)&\ge& C_{13}(J^a(u, v)\wedge 1).
\end{eqnarray*}
\end{lemma}

\pf
If $|u-v|\le r/2$, by Lemma \ref{lbbyskbm}
\begin{eqnarray*}
p^{a}_{B(u,r)\cup B(v,r)}(1/3, u, v) \ge \inf_{ |z|<r/2}
p^{a}_{B(0,r)}(1/3, 0, z) \ge c_1\left(r\wedge 1 \right)^2 e^{-c_2
r^2} \ge c_3\ge c_3(J^a(u, v)\wedge 1).
\end{eqnarray*}
Let $E= B(u,r)\cup B(v,r)$. If $|u-v|\ge r/2$, with  $E_1= B(u,r/8)$
and $E_3=B(v,r/8)$, we have by the strong Markov property and the
L\'evy system \eqref{e:levy} of $X^a$ that
\begin{eqnarray*}
p^{a}_{E}(1/3, u, v) &\ge&
\E_u
\left[p^{a}_{E}(1/3-\tau^{a}_{E_1},
X^{a}_{\tau^{a}_{E_1}}, v):\tau^{a}_{E_1}<1/3,
X^{a}_{\tau^{a}_{E_1}}\in E_3 \right]\\
&=& \int_0^{1/3} \left(\int_{E_1} p^a_{E_1}(s, u, w)
\left(\int_{E_3} J^a(w, z) p^a_E (1/3-s, z, v) dz\right) dw \right) ds \\
&\ge&\left(\inf_{w\in E_1,\, z\in E_3}J^{a}(w,z) \right)
\int_0^{1/3} \P_u \left(\tau^{a}_{E_1}>s \right) \left(\int_{E_3}
p^{a}_{E}(1/3-s, z, v) dz \right) ds\\
&\ge&  \P_u(\tau^{a}_{E_1}>1/3) \left(\inf_{w\in E_1,\, z\in
E_3}J^{a}(w,z) \right)\int_0^{1/3}\int_{E_3}p^{a}_{E_3}(1/3-s, z, v)dz ds\\
&=& \P_u(\tau^{a}_{E_1}>1/3)\left( \inf_{w\in E_1,\, z\in E_3}
j^{a}(|w-z|)\right)
\int_0^{1/3}\P_v(\tau^{a}_{E_3} > s)    ds\\
&\ge& \frac13 \P_u(\tau^{a}_{E_1}>1/3)\left( \inf_{w\in E_1,\, z\in
E_3} j^{a}(|w-z|)\right) \P_v(\tau^{a}_{E_3} > 1/3) \,.
\end{eqnarray*}
Thus by Lemma \ref{L:2.2},
\begin{eqnarray*}
p^{a}_{B(u,r)\cup B(v,r)}(1/3, u, v)&\ge&
\frac13\left(\P_0(\tau^{a}_{B(0,r/8)}>1/3)\right)^2 \left(\inf_{
w\in E_1,\,
z\in E_3}
j^{a}(|w-z|)\right)\\
&\ge& c_4 j^a(|u- v|) \ge c_4 (J^a(u, v) \wedge 1)\,.
\end{eqnarray*}
\qed

Recall that the function
$h^a_C(t,x,y)$ is defined in \eqref{eq:qd}.

\begin{thm}\label{t:low}
Suppose that $M$ and $T$ are positive constants. There are positive constants
$C_i=C_i(M, R_1, r_0, \lambda_0, \lambda, T, \alpha)$, $i=14, 15$,
such that for all $a\in (0, M]$ and $(t, x, y)\in (0, T] \times
D\times D$
\begin{equation}\label{eq:ppl2}
p^a_D(t, x, y)\ge C_{14}h^a_{C_{15}}(t,x,y).
\end{equation}
\end{thm}

\pf
Since $t^{-1/2}D$ satisfies the uniform interior ball condition
with radius $R_1(T)^{-1/2}$ for every $0<t \le T$,   there exist
$\delta=\delta(R_1, T) \in (0, R_1(T)^{-1/2})$ and  $L=L(R_1, T)>1$ such
that for all $t \in (0,  T]$ and $x,y \in D$, we can choose
$\xi_x \in (t^{-1/2}D)\cap B(t^{-1/2}x, L\delta)$ and $\xi_y \in
(t^{-1/2}D) \cap B(t^{-1/2}y, L\delta)$ with $ B(\xi_x, 2\delta)
\cap B(\xi_y, 2\delta)= \emptyset$ and $ B(\xi_x, 2\delta) \cup
B(\xi_y, 2\delta) \subset t^{-1/2}D$.

Let $x_t:=t^{-1/2}x$ and $y_t:=t^{-1/2}y$. Note that by
\eqref{e:direct} with $b=1/3$,
\begin{eqnarray}
\int_{B(\xi_x, \delta)} p^{a_t}_{{t^{-1/2}D}}(1/3,x_t,u)du &\ge& c_1
\left(\frac{\delta_D(x)}{\sqrt{t}} \wedge 1 \right) \int_{B(\xi_x,
\delta)} \left(\delta_{t^{-1/2}D}(u) \wedge 1 \right)
e^{-c_2|x_t-u|^2} du \nonumber\\
&\ge & c_1 \left(\frac{\delta_D(x)}{\sqrt{t}} \wedge 1 \right)
\, e^{-c_2 (L+1)^2 \delta^2}\,  |B(\xi_x, \delta)|
\nonumber \\
&\ge& c_3 \left(\frac{\delta_D(x)}{\sqrt{t}} \wedge 1 \right).
\label{e:loww_0}
\end{eqnarray}
Similarly
\begin{eqnarray}
\int_{B(\xi_y, \delta)} p^{a_t}_{{t^{-1/2}D}}(1/3,y_t,u)du \ge c_3
\left(\frac{\delta_D(y)}{\sqrt{t}} \wedge 1 \right).
\label{e:loww_01}
\end{eqnarray}

Now we deal with the cases $|x_t-y_t| \ge \delta/8$ and $|x_t-y_t| <
\delta/8$ separately.
Recall the definition of $a_t$ from \eqref{e:2.2}.

\medskip

\noindent
{\it Case 1}: Suppose $|x_t-y_t| \ge \delta/8$. Note that by the semigroup
property  and Lemma~\ref{lower bound12},
\begin{align*}
&p^{a_t}_{t^{-1/2}D}(1,x_t,y_t)\nonumber\\
\geq& \int_{B(\xi_y, \delta)}\int_{B(\xi_x, \delta)}
p^{a_t}_{t^{-1/2}D}(1/3,x_t,u)
p^{a_t}_{t^{-1/2}D}(1/3,u,v)p^{a_t}_{t^{-1/2}D}
(1/3,v,y_t)dudv \nonumber\\
\geq& \int_{B(\xi_y, \delta)}\int_{B(\xi_x, \delta)}
p^{a_t}_{{t^{-1/2}D}}(1/3,x_t,u)p^{a_t}_{B(u, \delta/2) \cup
B(v,\delta/2)}(1/3,u,v)p^{a_t}_
{t^{-1/2}D}(1/3,v,y_t)dudv\nonumber\\
\geq& c_4\int_{B(\xi_y, \delta)}\int_{B(\xi_x, \delta)}
p^{a_t}_{{t^{-1/2}D}}(1/3,x_t,u)(J^{a_t}(u,v)\wedge
1)p^{a_t}_{t^{-1/2}D}(1/3,v,y_t)dudv\nonumber\\
\geq& c_5 \left(\inf_{(u,v) \in B(\xi_x, \delta) \times
B(\xi_y,\delta)} (J^{a_t}(u,v)\wedge 1) \right)
\int_{B(\xi_y, \delta)}\int_{B(\xi_x, \delta)}
p^{a_t}_{{t^{-1/2}D}}(1/3,x_t,u)p^{a_t}_{t^{-1/2}D}(1/3,v,y_t)dudv.
\end{align*}
It then follows from \eqref{e:loww_0}--\eqref{e:loww_01} that
\begin{equation}\label{e:loww1}
p^{a_t}_{t^{-1/2}D}(1,x_t,
y_t) \ge  c_6 \left(\inf_{(u,v) \in
B(\xi_x, \delta) \times B(\xi_y,\delta)} (J^{a_t}(u,v)\wedge 1) \right)
\left(\frac{\delta_D(x)}{\sqrt{t}} \wedge 1 \right)
\left(\frac{\delta_D(y)}{\sqrt{t}} \wedge 1 \right).
\end{equation}
Using the fact that
\begin{equation}\label{e:jjj}
J^{a_t}(x_t, y_t) =a^{\alpha} t^{1+d/2}{\cal A}(d,
\alpha)|x-y|^{-(d+\alpha)}=t^{1+d/2}J^{a}(x,y)
\end{equation}
and the assumption $|x_t-y_t| \ge \delta/8$  which implies that
$|u-v| \le 2(1+L)\delta +|x_t-y_t| \le (17 +16L)|x_t-y_t|$,
we have
\begin{equation}\label{e:loww2}
\inf_{(u,v) \in B(\xi_x, \delta) \times B(\xi_y,\delta)}
(J^{a_t}(u,v)\wedge 1)\ge c_7\, (J^{a_t}(x_t,y_t)\wedge 1)=
c_7 \, ( t^{1+d/2} J^a(x,y)\wedge 1).
\end{equation}
Thus combining \eqref{e:loww1} and \eqref{e:loww2} with
\eqref{e:scaling},  we conclude that for $|x_t-y_t| \ge \delta/8$
\begin{eqnarray}
p^a_D(t, x, y)&=&t^{-d/2}p^{a_t}_{t^{-1/2}D}(1, t^{-1/2}x, t^{-1/2}y) \nonumber\\
&\ge& c_8t^{-d/2}\left(\frac{\delta_D(x)}{\sqrt{t}} \wedge 1
\right)\left(\frac{\delta_D(y)}{\sqrt{t}} \wedge 1 \right)
( t^{1+d/2} J^a(x,y)\wedge 1)\nonumber\\
&=& c_8\left(\frac{\delta_D(x)}{\sqrt{t}} \wedge 1 \right)
\left(\frac{\delta_D(y)}{\sqrt{t}} \wedge 1 \right) (
tJ^{a}(x,y)\wedge t^{-d/2}). \label{e:w2}
\end{eqnarray}

\noindent
{\it Case 2}: Suppose $|x_t-y_t| <\delta/8 $.
By the semigroup property,
\begin{align}
p^{a_t}_{t^{-1/2}D}(1,x_t,y_t) \geq  \int_{B(\xi_y,
\delta)}\int_{B(\xi_x, \delta)} p^{a_t}_{t^{-1/2}D}(1/3,x_t,u)
p^{a_t}_{t^{-1/2}D}(1/3,u,v)p^{a_t}_{t^{-1/2}D}(1/3,v,y_t)dudv
\label{e:sg11}.
\end{align}
By \eqref{e:direct} with $b=1/3$, we have for every $(u,v) \in
B(\xi_y, \delta) \times B(\xi_x, \delta)$,
$$
p^{a_t}_{t^{-1/2}D}(1/3,u, v) \ge c_9\left(\delta_{t^{-1/2}D}(u)
\wedge 1 \right)\left (\delta_{t^{-1/2}D}(v)\wedge 1 \right)
e^{-c_{10}|u-v|^2} \ge c_{11} (\delta\wedge 1)^2.
$$
Thus by  \eqref{e:loww_0}-\eqref{e:loww_01}
and \eqref{e:sg11},
\begin{eqnarray}
p^{a_t}_{t^{-1/2}D}(1,x_t,
y_t) &\ge& c_{12}
\left(\frac{\delta_D(x)}{\sqrt{t}}
\wedge 1 \right)\left(\frac{\delta_D(y)}{\sqrt{t}} \wedge 1 \right)\nonumber\\
&\ge& c_{12} \left(\frac{\delta_D(x)}{\sqrt{t}} \wedge 1 \right)
\left(\frac{\delta_D(y)}{\sqrt{t}} \wedge 1 \right) ( t^{1+d/2}
J^a(x,y)\wedge 1).\label{e:w1}
\end{eqnarray}

Combining \eqref{e:w2} and \eqref{e:w1} with Lemma \ref{lbbyskbm},
we have proved the theorem.      \qed

\section{Upper bound estimate} \label{s:ub}

In this section, we will establish upper bound estimate for
$X^a$ in any open set $D$ (not necessarily connected) satisfying
a weaker version of the uniform exterior ball condition.

Suppose that $U$ is a $C^{1,1}$ open set  with $C^{1,1}$
characteristics $(R, \Lambda)$. Without loss of generality, we can
always assume that $R\leq 1$ and $\Lambda \geq 1$. By definition,
for every  $Q\in \partial U$, there is a $C^{1,1}$-function $\phi_Q:
\bR^{d-1}\to \bR$ satisfying $\phi_Q(0)=0$, $ \nabla\phi_Q (0)=(0,
\dots, 0)$, $\| \nabla \phi_Q \|_\infty \leq \Lambda$, $| \nabla
\phi_Q (x)-\nabla \phi_Q (z)| \leq \Lambda |x-z|$, and an
orthonormal coordinate system $CS_Q$ $y=(\wt y, \, y_d)$ with origin
at $Q$ such that $B(Q, R )\cap U=\{ y=(\wt y, \, y_d)\in B(0, R)
\mbox{ in } CS_Q: y_d > \phi (\wt y) \}$. Define
$$
\rho_Q (x) := x_d -  \phi_Q (\wt x),
$$
where $(\wt x, x_d)$ is the coordinates of $x$ in $CS_Q$. Note that
for every $Q \in \partial U$ and $ x \in B(Q, R)\cap U$, we have
$(1+\Lambda^2)^{-1/2} \rho_Q (x) \le \delta_U(x) \le \rho_Q(x).$ We
define for $r_1, r_2>0$
$$
U_Q( r_1, r_2) :=\left\{ y\in U: r_1 >\rho_Q(y) >0,\, |\wt y | < r_2
\right\}.
$$

We recall the following key estimates from \cite[Lemma 3.5]{CKSV1}.

\begin{lemma}\label{L:2.0_1}
Suppose $R\in (0, 1]$, $M\in (0, \infty)$ and $\Lambda \in [1,
\infty)$ are constants, and let $r_0:=R/(4\sqrt{1+\Lambda^2})$.
There are constants $\delta_0 = \delta_0( R, M, \Lambda, \alpha)\in
(0, r_0)$, $C_{16}=C_{16} (R, M,\Lambda, \alpha)>0$
such that
for all $a \in (0, M]$,  $ \lambda\ge 1$, $C^{1,1}$ open set $U$
with characteristics $(R, \Lambda)$, $Q \in \partial U$ and $x \in
U_Q( \lambda^{-1} \delta_0 , \lambda^{-1} r_0 )$ with $\wt x =0$,
\bee\label{e:L:2.0_2}
\P_{x}\left(X^{a}_{\tau^{a}_{ U_Q( \lambda^{-1} \delta_0,
\lambda^{-1} r_0)}} \in  U\right)\le {C_{16}}{\lambda} \delta_U (x)
\eee
 and
\bee\label{e:L:2.0_3} \E_x\left[\tau^a_{  U_Q(
\lambda^{-1} \delta_0, \lambda^{-1} r_0)}\right]\,\le\,
{C_{16}}{\lambda^{-1}} \delta_U (x).
\eee
\end{lemma}

We note that
\begin{eqnarray*}
\P_x(\tau^{a}_{U}>1/4)& \leq&
\P_x \left(\tau^{a}_{ U_Q( \delta_0 , r_0)}>1/4 \right)+\P_{x}
\left(X^{a}_{\tau^{a}_{ U_Q( \delta_0,  r_0)}}
\in  U \hbox{ and } \tau^{a}_{ U_Q( \delta_0,  r_0)} \leq 1/4\right)\\
&\le &4 \, \E_x\left[\tau^a_{  U_Q( \delta_0,  r_0)}\right]+
\P_{x}\left(X^{a}_{\tau^{a}_{ U_Q( \delta_0,  r_0)}} \in  U\right).
\end{eqnarray*}
Thus, by \eqref{e:L:2.0_2}-(\ref{e:L:2.0_3}) with $\lambda=1$ and a simple geometric consideration,
we obtain the following lemma.

\begin{lemma}\label{lem:etd1r_2}
Suppose that $M> 0$ and $U$ is a $C^{1, 1}$  open set with the
characteristics $(R, \Lambda)$.
There exists $C_{17}=C_{17}
(\Lambda, R, M, \alpha)>0$
such that for all
$a\in (0, M]$ and $x\in U$,
$$
\P_x(\tau^{a}_{U}>1/4) \le C_{17}\delta_U(x)
.
$$
\end{lemma}

In particular, we have the following.

\begin{cor}\label{C:3.3}
Suppose that $M$ and $r_1$ are positive constants and $E := \{x\in
\bR^d: \, |x-x_0| >r_1\}$.
There exists $C_{18}=C_{18}
(r_1, M, \alpha)>0$
independent of $x_0$ such
that for all $a\in (0, M]$ and $x\in E$,
$$
\P_x(\tau^{a}_{E}>1/4) \le C_{18}\delta_E(x)
.
$$
\end{cor}

The proof of the next lemma is  similar to that of \cite[Lemma
2]{BGR}, which is a variation of the proof of \cite[Lemma 2.2]{CKS}.
We give the proof here for the sake of completeness.

\begin{lemma}\label{l:gen}
Suppose that $E_1,E_3, E$ are open subsets of $\bR^d$ with $E_1,
E_3\subset E$ and ${\rm dist}(E_1,E_3)>0$.
For any $n\ge 1$, let $E_{2,i}$, $i=1, \dots, n$, be disjoint Borel subsets with
$\cup_{i=1}^nE_{2,i} =E\setminus
(E_1\cup E_3)$.  If $x\in E_1$ and $y \in E_3$, then for
all $a > 0$ and $t >0$,
\begin{equation}\label{eq:ub}
p^{a}_{E}(t, x, y) \le \sum_{i=1}^n\P_x(X^{a}_{\tau^{a}_{E_1}}\in E_{2,i})
\left(\sup_{s<t,\, z\in E_{2,i}}
p_E^a(s, z, y)\right)+
(t \wedge \E_x [\tau^{a}_{E_1}])
\left(\sup_{u\in E_1,\, z\in E_3}J^{a}(u,z)\right).
\end{equation}
\end{lemma}

\pf Using the strong Markov property, we have
\begin{eqnarray*}
p^{a}_{E}(t, x, y) &=&\E_x\left[p^{a}_{E}\big(t-\tau^{a}_{E_1},
X^{a}_{\tau^{a}_{E_1}}, y\big)
: \tau^{a}_{E_1}<t \right]\\
&=&\sum_{i=1}^n \E_x\left[p^{a}_{E}\big(t-\tau^{a}_{E_1},
X^{a}_{\tau^{a}_{E_1}},
y\big):\tau^{a}_{E_1}<t, X^{a}_{\tau^{a}_{E_1}}\in E_{2,i} \right] \\
&&~~+ \E_x\left[p^{a}_{E}\big(t-\tau^{a}_{E_1},
X^{a}_{\tau^{a}_{E_1}}, y\big):\tau^{a}_{E_1}<t,
X^{a}_{\tau^{a}_{E_1}}\in E_3\right] \,=:\, I\,+\,II\,.
\end{eqnarray*}
Clearly
\begin{eqnarray*}
I &\le&  \sum_{i=1}^n \P_x\left(\tau^{a}_{E_1}<t, X^{a}_{\tau^{a}_{E_1}}\in
E_{2,i}\right) \left( \sup_{s<t,\, z\in E_{2,i}} p_E^a(s, z, y)\right) \\
&\le& \sum_{i=1}^n \P_x\left( X^{a}_{\tau^{a}_{E_1}}\in
E_{2,i}\right) \left( \sup_{s<t,\, z\in E_{2,i}} p_E^a(s, z,
y)\right) .
\end{eqnarray*}
On the other hand,
by \eqref{e:levy},
\begin{eqnarray*}
II&=& \int_0^{t} \left( \int_{E_1} p^{a}_{{E_1}}(s, x, u) \left(
\int_{E_3} J^{a}(u,z) p^{a}_{E}(t-s, z, y) dz\right)
du\right) ds\\
&\le&   \left(\sup_{u\in E_1,\, z\in E_3}J^{a}(u,z)\right)
\int_0^{t} \P_x(\tau^{a}_{E_1}>s) \left(\int_{E_3}p^{a}_{E}(t-s, z,
y)dz
\right) ds\\
&\le& \int_0^{t}\P_x(\tau^{a}_{E_1}>s) ds \sup_{u\in E_1,\, z\in
E_3}J^{a}(u,z) \le
(t \wedge \E_x [\tau^{a}_{E_1}])
\sup_{u\in E_1,\, z\in
E_3}J^{a}(u,z)\,.
\end{eqnarray*}
This completes the proof of the lemma. \qed

\begin{thm}\label{ub11}
Suppose that $M>0$ is a constant and that
$D$ is an open set satisfying a weaker
version of the uniform exterior ball condition with radius $R_1>0$.
There exists a positive constant $C_{19}=C_{19}(M, R_1)$ such that for
all $a \in (0, M]$ and $x, y\in D$,
\begin{equation}\label{eq:ppu}
p^{a}_{D}(1/2,x,y)\leq C_{19}\,
(\delta_D(x) \wedge 1)\,
\left(e^{-|x-y|^2/(2C_2)} +\left(
j^a(|x-y|)\wedge 1\right) \right)\, .
\end{equation}
\end{thm}

\pf
First note that for every $x_0\in \bR^d$, $\{z \in \bR^d:
|x-x_0|> R_1 \}$ is
a $C^{1,1}$ open set with characteristics
$(R, \Lambda)$ depending only on $R_1$ and $d$. Let $r_0$ and $\delta_0$
be the positive constants in Lemma \ref{L:2.0_1} for $U=\{z \in
\bR^d: |x-x_0|> R_1\}$.

It follows from Corollary \ref{C:1.2} that
$$
p^a_{D} (1/2, x, y) \, \leq   \, p^a (1/2, x, y)\,\le
c_1\left(e^{-|x-y|^2/(2C_2)}+
\left(j^a(|x-y|)\wedge 1\right)
\right),
$$
so it suffices to prove the theorem for $x\in D$  with
$\delta_{D} (x) <\delta_0/(32)$.

Now fix $x \in D$ with $\delta_{D} (x) <\delta_0/(32)$ and let $Q\in
\partial D$ be such that $|x-Q|=\delta_{D} (x)$. Let $B_Q\subset
D^c$ be the ball with radius
$R_1$ so that
$\partial B_Q \cap \partial D=\{Q\}$
and $E:=(\overline{B_Q})^c$.
Observe that $\delta_E(x)=\delta_D(x)=|x-Q|$.

When
$|x-y|\le \sqrt{dC_2} \vee
((\delta_0+r_0)/2)$, we have from  Corollary \ref{C:1.2} that
$$
p^a(1/2,x,y) \ge c_2e^{-c_3|x-y|^2}  \ge c_4 >0 \quad
\text{and}\quad \sup_{z\in \bR^d} p^a(1/4,z,y) \le c_5.
$$
Thus, by the semigroup property  and Corollary \ref{C:3.3},
\begin{eqnarray}
p^{a}_{D}(1/2,x,y)&=&
\int_D p^a_D(1/4,x,z)p^a_D(1/4,z,y)dz \nonumber\\
& \leq&
\sup_{z\in D} p_D^a(1/4,z,y) \P_x(\tau^{a}_{D}>1/4)\nonumber\\
& \leq& \sup_{z\in \bR^d} p^a(1/4,z,y) \P_x(\tau^{a}_{E}>1/4)\nonumber\\
&\leq& c_6 \delta_{E}(x)=c_6 \delta_{D}(x)\le
c_7\,\delta_{D}(x)p^a(1/2,x,y)\,.\label{eq:ppun2}
\end{eqnarray}

Finally we consider the case that
$|x-y|>\sqrt{dC_2} \vee ((\delta_0+r_0)/2)$
(and $\delta_{D} (x) <\delta_0/(32)$).

There is a $C^{1,1}$-function $\phi: \bR^{d-1}\to \bR$ satisfying
$\phi(0)=0$, $ \nabla\phi (0)=(0, \dots, 0)$,
$\| \nabla \phi  \|_\infty \leq
\Lambda$, $| \nabla \phi (w)-\nabla \phi (z)| \leq \Lambda |w-z|$,
and an orthonormal coordinate system $CS$ with its origin at $Q$
such that
$$
B(Q, R)\cap E=\{ z=(\wt z, z_d) \in B(0, R) \mbox{ in } CS: z_d >
\phi (\wt z) \}
$$
and that $x$ has coordinate $(\wt 0, \delta_D(x))$ in  this $CS$.
Let
$$
E_1:=\left\{ z=(\wt z, z_d) \mbox{ in } CS: 0<z_d-\phi(\wt
z)<\delta_0/8, \ |\wt z| < r_0/8\right\},
$$
$E_3:=
\{z\in E: |z-x|>|x-y|/2\}$
and $E_2:=E\setminus (E_1\cup E_3)$. Note that
$|z-x| >(\delta_0+r_0)/4 $ for $z\in E_3$. So,
if $u\in E_1$ and $z\in E_3$,
then
\begin{eqnarray}\label{e:n00}
|u-z| \ge |z-x|-|x-u| \ge |z-x|-(\delta_0+r_0)/8 \ge
\frac{1}{2}|z-x| \ge \frac{1}{4}|x-y|.
\end{eqnarray}
Thus
\begin{eqnarray}
\sup_{u\in E_1,\, z\in E_3}J^{a}(u,z)
&\le& \sup_{(u,z):|u-z| \ge
\frac{1}{4}|x-y|}J^{a}(u,z)  \nonumber \\
&\le&
j^{a}(|x-y|/4) = \left(j^{a}(|x-y|/4)\wedge  j^M(    (\delta_0+r_0)/8   )
\right).   \label{e:n0}
\end{eqnarray}
If $z \in E_2$, then $|z-y| \ge |x-y| -|x-z| \ge |x-y|/2$.
We also observe  that for every $\beta \ge d/4$,
$\sup_{s <1/2} s^{-d/2} e^{-\beta/s} =2^{d/2}e^{-2\beta}$.
By Corollary \ref{C:1.2} and these observations,
\begin{eqnarray}
\sup_{s<1/2,\, z\in E_2} p^a(s, z, y)&\le& C_3  \sup_{s<1/2,\, z \in
E_2}\left(s^{-d/2}e^{-|z-y|^2/(C_2s)}+ \big(s^{-d/2}
\wedge s J^a(z,y)\big)\right)
\nonumber\\
&=& C_3 2^{d/2} e^{-|x-y|^2/C_2}+ \frac{C_3}2   j^a(|x-y|/2)
\nonumber\\
&=&
C_3 2^{d/2}e^{-|x-y|^2/C_2} +\frac{C_3}2 \left(
j^a(|x-y|/2) \wedge  j^M(    (\delta_0+r_0)/4   )    \right)\nonumber\\
&\le&  c_8 \left(e^{-|x-y|^2/C_2} +  \big( j^a(|x-y|) \wedge 1\big)
\right) \label{e:n1}
\end{eqnarray}
for some $c_8>0$. Applying Lemmas \ref{L:2.0_1} and \ref{l:gen}, we
obtain,
\begin{eqnarray*}
p^{a}_{E}(1/2, x, y)
&\le& c_9 \left(e^{-|x-y|^2/(2C_2)} +\big(
j^a(|x-y|)\wedge 1\big)\right)\left(\P_x(X^{a}_{\tau^{a}_{E_1}}\in
E)+
\E_x[\tau^{a}_{E_1}]\right)\\
&\leq& c_{10}\, \delta_E(x)\left(e^{-|x-y|^2/(2C_2)} +
\big( j^a(|x-y|)\wedge 1 \big) \right)\\
&=& c_{10} \,\delta_D(x)\left(e^{-|x-y|^2/(2C_2)} +
\big( j^a(|x-y|)\wedge 1 \big) \right).
\end{eqnarray*}
Therefore
$$
p^{a}_{D}(1/2, x, y) \le p^{a}_{E}(1/2, x, y)\le c_{10}
\delta_D(x)\left(e^{-|x-y|^2/(2C_2)} +\big( j^a(|x-y|)\wedge
1\big)\right).
$$
\qed

\begin{thm}\label{t:ub}
Assume that $M>0$ is a constant and
that $D$ is an open set satisfying a weaker
version of the uniform exterior ball condition with radius $R_1>0$.
For every $T>0$, there exists a positive constant $C_{20}=C_{20}(T,
R_1, \alpha, M)$ such that for all $a \in [0, M]$
and $(t, x, y)\in (0, T]\times D \times D$,
\bee\label{e:1} p^a_D(t, x, y) \leq C_{20} \left( 1\wedge
\frac{\delta_D (x)}{\sqrt{t}}\right) \left( 1\wedge
\frac{\delta_D(y)}{\sqrt{t}}\right)
\left( t^{-d/2}
e^{-|x-y|^2/(4C^3_2t)}+ \big(t^{-d/2} \wedge tJ^a(x, y)\big) \right).
\eee
\end{thm}

\pf Fix $T, M>0$ and recall that
$a_t:=
a t^{(2-\alpha)/(2\alpha)} \le
M T^{(2-\alpha)/(2\alpha)}$. Note that $t^{-1/2}D$ is an open set
satisfying a weaker version of the uniform exterior ball condition
with radius $T^{-1/2}R_1>0$ for every $t \in (0, T]$. Thus, by
Theorem \ref{ub11}, there exists a positive constant $c_1=c_1(T,
R_1, \alpha, M)$ such that for all $ t \in (0, T]$ and $a \in (0,
M]$,
\bee\label{e:4}
p^{a_t }_{t^{-1/2}D}
(1/2,x,y)
\leq c_1  (e^{-|x-y|^2/(2C_2)} +(
j^{a_t}(|x-y|)\wedge 1))
\delta_{t^{-1/2}D}(x).
\eee
Thus
by \eqref{e:scaling}, \eqref{e:jjj} and \eqref{e:4}, for every
$t \le T$,
\begin{eqnarray*}
p^a_D(t/2, x, y) &=& t^{-d/2} p^{a_t}_{t^{-1/2}D} (1/2, t^{-1/2}x, t^{-1/2} y)\\
&\leq & c_1 \, t^{-d/2}  \left(e^{-|x-y|^2/(2C_2t)} +(
j^{a_t}(|x-y|/t^{1/2})\wedge 1)\right) \delta_{t^{-1/2} D}(t^{-1/2}
x)\\
&=& c_1  \left( t^{-d/2}  e^{-|x-y|^2/(2C_2t)}+ (t^{-d/2}
\wedge tJ^a(x, y)) \right)
\frac{\delta_{D}( x)}{\sqrt{t}} .
\end{eqnarray*}
By  symmetry, the  above  inequality  holds with the roles of $x$
and $y$ interchanged.
Using the semigroup property
and Corollary \ref{C:1.2} (twice),
for $t\leq T$,
\begin{eqnarray*}
p^a_D(t, x, y) &=& \int_D p^a_D(t/2, x, z) p^a_D (t/2, z, y) dz\\
&\leq & c_3\,  \frac{\delta_{D}( x)\delta_{D}(
y)}{t } \int_D p^a(2C^2_2t , x, z) p^a (2C^2_2 t, z, y) dz \\
&\leq & c_3\,  \frac{\delta_{D}( x)\delta_{D}(
y)}{t }  p^a(4C^2_2t , x, y)\\
&\leq & c _4\, \frac{\delta_{D}( x) \delta_{D}( y)}{t } \left(
t^{-d/2}  e^{-|x-y|^2/(4C^3_2t)}+ (t^{-d/2} \wedge tJ^a(x, y))
\right).
\end{eqnarray*}
This with Corollary \ref{C:1.2} proves the upper bound \eqref{e:1}
by noting that
$$
(1\wedge u)(1\wedge v) = \min\{1, u, v, uv\} \qquad \hbox{for } u,
v>0.
$$
\qed

We point out that,
in view of Theorem \ref{t:low}, the above upper bound estimate \eqref{e:1}
is sharp when $x$ and $y$ are in the same component of $D$.
However it is not sharp when $x$ and $y$ are
in
different components of $D$, since in this case
when $a\to 0$, it does not go to zero  and thus does not
give the sharp upper bound
for the Dirichlet heat kernel $p_D^0(t, x, y)$ of Brownian motion in $D$.
Next we improve the above estimate to get the sharp estimate stated in
Theorem \ref{t:ub1} below.

For the remainder of this section, we continue assume $D$ is an open set satisfying a weaker
version of the uniform exterior ball condition with radius $R_1>0$.
It is easy to see that
the distance between any two distinct
connected components of $D$ is at least $R^*$ for some $R^*>0$
that depends only on $R_1$. Without
loss of generality, we assume that $R^*=R_1$.
Observe that for
$c_0>0$, $r\geq r_0$ and $t>0$,
\begin{equation}\label{e:3.11}
t^{-d/2}e^{-c_0 r^2/t} \leq c_1 t^{-d/2} (t/r^2)^{d/2+1} = c_1 \frac{t}{r^{d+2}}
\leq c_1 r_0^{\alpha -2}   \,  \frac{t}{r^{d+\alpha}},
\end{equation}
where $c_1>0$ depends only on $c_0$, $r_0$ and
and $d$.
This implies that for $x$ and $y$ in different components of $D$,
the jumping kernel component
$tJ^1(x, y)$ dominates the Gaussian component $t^{-d/2} e^{-|x-y|^2/C_2t}$.
This fact will be used several times in the rest of this section.

By Theorem \ref{t:ub}, we only need to consider the case when $x$
and $y$ are in different components of $D$.
Recall that , for any $x\in D$, $D(x)$ denotes the connected
component  of $D$ that contains $x$.

First we give an interior upper bound of $p^{a}_{D}(t, x, y)$ when
$x$ and $y$ are in different components of $D$.

\begin{lemma}\label{l:gen1}
Assume that $M>0$ is a constant and
that $D$ is an open set satisfying a weaker
version of the uniform exterior ball condition with radius $R_1>0$.
For every $T>0$, there exists a positive constant $C_{21}=C_{21}(T,
R_1, \alpha,   M)$ such that for
all $a \in (0, M]$, $t \in (0, T]$ and $x, y$ in different components of $D$,
$$p^{a}_{D}(t, x, y) \le C_{21}
a^\alpha t |x-y|^{-d-\alpha}.
$$
\end{lemma}

\pf Using the strong Markov property and \eqref{e:levy}, we have for $t \le T$,
\begin{eqnarray}
&&p^{a}_{D}(t, x, y) \nonumber\\
&=&
\E_x\left[p^{a}_{D}\big(t-\tau^{a}_{D(x)},
X^{a}_{\tau^{a}_{D(x)}}, y\big):\tau^{a}_{D(x)}<t,
X^{a}_{\tau^{a}_{D(x)}}\in D\setminus D(x)\right]\nonumber\\
&=& \int_0^{t} \left( \int_{D(x)} p^{a}_{{D(x)}}(s, x, u) \left(
\int_{D\setminus D(x)} J^{a}(u,z) p^{a}_{D}(t-s, z, y) dz\right)
du\right) ds\nonumber\\
&\le&
c_1\frac{a^\alpha }{t M^\alpha} \int_0^{t} \left( \int_{D(x)} p^{a}(s, x, u) \left(
\int_{D\setminus D(x)} \left(t^{-d/2} \wedge (tJ^{M}(u,z))\right) p^{a}(t-s, z, y) dz
\right) du\right) ds\nonumber\\
&\le& c_1\frac{a^\alpha }{tM^\alpha} \int_0^{t} \left( \int_{\bR^d}
p^{a}(s, x, u) \left( \int_{\bR^d} \left( t^{-d/2} \wedge (tJ^{M}(u,z))
\right) p^{a}(t-s, z, y) dz\right) du\right) ds . \label{e:ew1}
\end{eqnarray}
In the second to the last inequality above,
we have used the facts that
$J^{M}(u,z) \le j^M(R_1)$ and $t \le T$. By Corollary \ref{C:1.2},
$p^{a}(s, x, u) \le  c_2 p^{M} (C_2^2s,x,u)$, $p^{a}(t-s, z, y) \le
c_2 p^M (C_2^2(t-s), z,y)$ and $ t^{-d/2} \wedge (tJ^{M}(u,z)) \le
c_2 p^{M} (t,u,z).$
Thus, using the semigroup property and Corollary \ref{C:1.2}, from
\eqref{e:ew1} we obtain that
\begin{eqnarray*}
p^{a}_{D}(t, x, y)
&\le& c_3
\frac{a^\alpha }{t M^\alpha} \int_0^{t} \left( \int_{\bR^d} p^{M} (C_2^2s,x,u)\left(
\int_{\bR^d} p^{M} (t,u,z) p^M (C_2^2(t-s), z,y) dz\right)
du\right) ds\\
&=&  c_3 \frac{a^\alpha }{t M^\alpha}
\int_0^t  p^M ((C_2^2+1)t, x,y) ds \\
&\le&
c_4  a^\alpha \left(  t^{-d/2}
e^{-|x-y|^2/(C_2(C_2^2+1)t)}+ t^{-d/2}\, \wedge
\frac{t}{|x-y|^{d+\alpha}}  \right)\\
&\le&   c_5 \,
a^\alpha \,t \,|x-y|^{-d-\alpha}.
\end{eqnarray*}
In the last inequality above,
we have used the fact that
$|x-y| \ge R_1$.
\qed

\begin{thm}\label{ub111}
Assume that $M>0$ is a constant and
that $D$ is an open set satisfying a weaker
version of the uniform exterior ball condition with radius $R_1>0$.
There exists a positive constant
$C_{22}=C_{22}(M, R_1 )$
 such that for all $a \in (0, M]$ and $x, y$ in different components of $D$.
$$
p^{a}_{D}(1 , x, y)
\leq C_{22}\,(\delta_D(x) \wedge 1)(\delta_D(y) \wedge 1)
\left(
j^a(|x-y|)\wedge 1\right)\, .
$$
\end{thm}

\pf
We first claim that
\begin{equation}\label{e:dfd1}
p^{a}_{D}(1/2, x, y) \le
c_{1}a^\alpha (\delta_D(x) \wedge 1)\left(
|x-y|^{-d-\alpha}\wedge 1\right).
\end{equation}
Recall  that for every $x_0\in \bR^d$, $\{z \in \bR^d:
|x-x_0|>  R_1/4 \}$
is a $C^{1,1}$ open set with characteristics
$(R, \Lambda)$ depending only on $R_1$ and $d$. Let $r_0$ and $\delta_0$
be the positive constants in Lemma \ref{L:2.0_1} for $U=\{z \in
\bR^d: |x-x_0|> R_1/4\}$.
It follows from Lemma \ref{l:gen1} that
$$
p^a_{D} (1/2, x, y)   \le
c_1 \left(j^a(|x-y|)\wedge 1\right).
$$
So it suffices to prove \eqref{e:dfd1} for $x\in D$  with
$\delta_{D} (x) <\delta_0/(32)$.

Now fix $x \in D$ with $\delta_{D} (x) <\delta_0/(32)$ and let $Q\in
\partial D$ be such that $|x-Q|=\delta_{D} (x)$. Let $B_Q:=B(
Q, R_1/4)\subset D^c$ be the ball with radius $R_1/4$ so that
$\partial B_Q \cap \partial D=\{Q\}$.

There is a $C^{1,1}$-function $\phi: \bR^{d-1}\to \bR$ satisfying
$\phi(0)=0$, $ \nabla\phi (0)=(0, \dots, 0)$,
$\| \nabla \phi  \|_\infty \leq
\Lambda$, $| \nabla \phi (w)-\nabla \phi (z)| \leq \Lambda |w-z|$,
and an orthonormal coordinate system $CS$ with its origin at $Q$
such that
$$
B(Q, R)\cap  (\overline{B_Q})^c=\{ z=(\wt z, z_d) \in B(0, R) \mbox{ in } CS: z_d >
\phi (\wt z) \}
$$
and that $x$ has coordinate $(\wt 0, \delta_D(x))$ in  this $CS$.
Let
$$
E:=D \cup ( B(
Q, R_1/2)  \setminus \overline{B_Q}),$$
$$
E_1:=\left\{ z=(\wt z, z_d) \mbox{ in } CS: 0<z_d-\phi(\wt
z)<\delta_0/8, \ |\wt z| < r_0/8\right\},
$$
$E_3:=
\{z\in E: |z-x|>|x-y|/2\}$, $E_{2,1}:=(E\setminus (E_1\cup E_3)) \cap D(y)$ and
$E_{2,2}:=  E\setminus (E_1\cup E_3  \cup D(y))$.
Observe that $\delta_{(\overline{B_Q})^c}(x)=\delta_{E}(x)=\delta_D(x)=|x-Q|$ and
$|x-y| \ge R_1 > R >
((\delta_0+r_0)/2)$.
So, by \eqref{e:n00}--\eqref{e:n0}, dist$(E_1, E_3) \ge R_1/4$ and
\begin{eqnarray}
\sup_{u\in E_1,\, z\in E_3}J^{a}(u,z) \le c_2 \left(j^{a}(|x-y|)\wedge 1   \right). \label{e:f1}
\end{eqnarray}

If $z \in E_{2,i}, i=1,2$, then $|z-y| \ge |x-y| -|x-z| \ge |x-y|/2$. Thus
by the same argument as the one in  \eqref{e:n1}, if
$|x-y|>\sqrt{dC_2}$,
we have by \eqref{e:3.11}
\begin{eqnarray}
\sup_{s<1/2,\, z\in E_{2,1}} p^a(s, z, y)\le   c_3 \left(e^{-|x-y|^2/C_2} +\big( j^a(|x-y|) \wedge 1\big)
\right)\le   c_4 \left( j^M(|x-y|) \wedge 1 \right).\label{e:f2}
\end{eqnarray}
If
$|x-y| \le\sqrt{dC_2} $, since $R_1 \le |x-y|$, we also have
\begin{eqnarray}
\sup_{s<1/2,\, z\in E_{2,1}} p^a(s, z, y)&\le& C_3 \sup_{s<1/2,\, z\in E_{2,1}}
\left(s^{-d/2}e^{-|z-y|^2/(C_2s)}+ \big(s^{-d/2}
\wedge s J^a(z,y)\big)\right)\nonumber\\
&\le & c_5\sup_{s<1/2}\left(s^{-d/2}e^{-R_1^2/(2C_2s)}+ \big(s^{-d/2}
\wedge s J^a(R_1)\big)\right)\nonumber\\
&\le & c_6 \,\le\, c_7\,  j^M(|x-y|).\label{e:f2-1}
\end{eqnarray}
On the other hand, since $D(y) \subset E_{2,2}^c $, by Lemma \ref{l:gen1},
\begin{eqnarray}
\sup_{s<1/2,\, z\in E_{2,2}} p_E^a(s, z, y)\le
C_{21} \sup_{s<1/2,\, z\in E_{2,2}}  s J^a(|z-y|)  \le  c_8 \big( j^a(|x-y|) \wedge 1\big).\label{e:f3}
\end{eqnarray}
Furthermore, since ${\rm dist} (E_1, D(y))\ge R_1/2$,
by the L\'evy system \eqref{e:levy},
\begin{eqnarray}
\P_x(X^{a}_{\tau^{a}_{E_1}}\in
{E_{2,1}}) &\le&
\P_x(X^{a}_{\tau^{a}_{E_1}}\in
{D(y)}) \nonumber\\
&=&
\int_0^{\infty} \left( \int_{E_1} p^{a}_{{E_1}}(s, x, u) \left(
\int_{D(y)} J^{a}(u,z)  dz\right)
du\right) ds\nonumber\\
&=& a^\alpha
\int_0^{\infty} \left( \int_{E_1} p^{a}_{{E_1}}(s, x, u) \left(
\int_{D(y)} J^{1}(u,z)  dz\right)
du\right) ds\nonumber\\
&\le & a^\alpha \left(
\int_{\{|z|>
R_1/2\}} J^{1}(z)  dz\right)
\int_0^{\infty} \left( \int_{E_1} p^{a}_{{E_1}}(s, x, u) du\right)ds \nonumber \\
&\le&  c_9 a^\alpha \E_x[\tau^{a}_{E_1}] .\label{e:f4}
\end{eqnarray}

Applying Lemmas \ref{L:2.0_1} and \ref{l:gen}, and combining \eqref{e:f1}--\eqref{e:f4},  we
obtain,
\begin{eqnarray*}
&&p^{a}_{E}(1/2, x, y)\\
&\le&
\P_x(X^{a}_{\tau^{a}_{E_1}}\in E_{2,1})
\left(\sup_{s<1/2,\, z\in E_{2,1}}
p_E^a(s, z, y)\right)+  \P_x(X^{a}_{\tau^{a}_{E_1}}\in E_{2,2})
\left(\sup_{s<1/2,\, z\in E_{2,2}}
p_E^a(s, z, y)\right)\\
&&\quad+
\E_x [\tau^{a}_{E_1}]
\left(\sup_{u\in E_1,\, z\in E_3}J^{a}(u,z)\right)\\
&\le&
c_9 a^\alpha \E_x [\tau^{a}_{E_1}]
\left(\sup_{s<1/2,\, z\in E_{2,1}}
p^a(s, z, y)\right)+  \P_x(X^{a}_{\tau^{a}_{E_1}}\in (\overline{B_Q})^c)
\left(\sup_{s<1/2,\, z\in E_{2,2}}
p_E^a(s, z, y)\right)\\
&&\quad+
\E_x [\tau^{a}_{E_1}]
\left(\sup_{u\in E_1,\, z\in E_3}J^{a}(u,z)\right)\\
&\leq& c_{10}\, a^\alpha(\delta_D(x) \wedge 1)\left(
|x-y|^{-d-\alpha}\wedge 1\right).
\end{eqnarray*}
Therefore we have proved the claim \eqref{e:dfd1}.
In particular, we have
\begin{equation}\label{e:dfd}
p^{a}_{D}(1/2, x, y) \le c_{10}
a^\alpha(\delta_D(x) \wedge 1)\left(e^{-|x-y|^2/(2C_2)} + \left(
|x-y|^{-d-\alpha}\wedge 1\right)\right).
\end{equation}
By  symmetry, the  above  inequality  holds with the roles of $x$
and $y$ interchanged.
It follows from the semigroup property that
\begin{eqnarray*}
&&p^a_D(1, x, y)\\
&=& \int_D p^a_D(1/2, x, z) p^a_D (1/2, z, y) dz\\
&=&\int_{D(x)} p^a_D(1/2, x, z) p^a_D (1/2, z, y) dz+ \int_{D \setminus D(x)}
p^a_D(1/2, x, z) p^a_D (1/2, z, y) dz.
\end{eqnarray*}
By applying Theorem \ref{ub11} and
\eqref{e:dfd}, and then  applying Corollary \ref{C:1.2} (twice),  we have
\begin{eqnarray*}
p^a_D(1, x, y)
&\leq & c_{11}\,a^{\alpha} (\delta_D(x) \wedge 1) (\delta_D(y) \wedge 1)\int_D p^1(2C_2^2, x, z) p^1 (2C_2^2, z, y) dz \\
&\leq & c_{11} \, a^{\alpha}(\delta_D(x) \wedge 1) (\delta_D(y) \wedge 1) p^1(4C_2^2,x,y)\\
&\leq & c_{12} \, a^{\alpha}(\delta_D(x) \wedge 1) (\delta_D(y) \wedge 1) \left(e^{-|x-y|^2/(4C^3_2)} + \left(
|x-y|^{-d-\alpha}\wedge 1\right)\right)\\
&\leq & c_{13} \, a^{\alpha}(\delta_D(x) \wedge 1) (\delta_D(y) \wedge 1)
|x-y|^{-d-\alpha}
\\
&\leq & c_{14} \, (\delta_D(x) \wedge 1) (\delta_D(y) \wedge 1)\left(
J^a(x, y) \wedge 1 \right).
\end{eqnarray*}
\qed

Recall that $h^a_C(t,x,y)$ is defined in \eqref{eq:qd}.

\begin{thm}\label{t:ub1}
Assume that $M>0$ is a constant and
that $D$ is an open set satisfying a weaker
version of the uniform exterior ball condition with radius $R_1>0$.
For every $T>0$, there exists a positive constant $C_{23}=C_{23}(T,
R_1, \alpha,  M)$ such that for all $a \in [0, M]$
and $(t,x,y)\in (0, T] \times D \times D$,
$$ p^a_D(t, x, y) \,\leq \,C_{23} \,
h^a_{ (2C_2)^{-1}}(t,x,y).
$$
\end{thm}

\pf Fix $T, M>0$ and recall that $a_t=a t^{(2-\alpha)/(2\alpha)} \le
M T^{(2-\alpha)/(2\alpha)}$. Note that $t^{-1/2}D$ is an open set
satisfying a weaker version of the uniform exterior ball condition
with radius $T^{-1/2}R_1>0$ for every $t \in (0, T]$. Thus, by
Theorem \ref{ub111}, there exists a positive constant $c_1=c_1(T,
R_1, \alpha,   M)$ such that for all $ t \in (0, T]$, $a \in (0,
M]$ and
$x, y$ in different components of $D$,
\bee\label{e:3.20}
p^{a_t }_{t^{-1/2}D} (1, x, y)
\leq c_1  \left( j^{a_t}(|x-y|)\wedge 1 \right)
\left(\delta_{t^{-1/2}D}(x) \wedge 1\right) \left(\delta_{t^{-1/2}D}(y) \wedge 1\right).
\eee
Thus
by \eqref{e:scaling}, \eqref{e:jjj} and
\eqref{e:3.20},
for every $t \le T$, $a \in (0, M]$ and
$x, y$ in different components of $D$,
\begin{eqnarray*}
p^a_D(t, x, y) &=& t^{-d/2} p^{a_t}_{t^{-1/2}D} (1, t^{-1/2}x,
t^{-1/2} y)\\
&\leq & c_1 \, t^{-d/2}  \left(
j^{a_t}(|x-y|/t^{1/2})\wedge 1\right) \left(\delta_{t^{-1/2} D}(t^{-1/2}
x)\wedge 1 \right) \left(\delta_{t^{-1/2} D}(t^{-1/2}
y)\wedge 1 \right)\\
&=& c_1  \left( t^{-d/2}
\wedge tJ^a(x, y) \right)
\left(\frac{\delta_{D}( x)}{\sqrt{t}}\wedge 1 \right)
\left( \frac{\delta_{D}( y)}{\sqrt{t}}\wedge 1 \right) .
\end{eqnarray*}
Combining this result
with Theorem \ref{t:ub}, we have proved the theorem.
\qed

\section{Large time estimates} \label{s:l}

In this section we assume that $D$ is a bounded $C^{1,1}$ open set
in $\bR^d$ which may be disconnected, and we give the proof of
Theorem \ref{t:main} (ii) and (iii).

\medskip

\noindent {\bf Proof of  Theorem \ref{t:main} (ii) and (iii).}
Let $D$ be a bounded $C^{1,1}$ open set in $\bR^d$ with $d\geq 1$.
For each $a\geq 0$, the semigroup of $X^{a, D}$ is Hilbert-Schmidt
as, by Theorem \ref{T:1.1}
$$\int_{D\times D} p^a_D (t, x, y)^2 dx dy= \int_D p^a_D(2t, x, x) dx
\le  C_1\, \,\left((2t)^{-d/2} \wedge
(a^{\alpha}2t)^{-d/\alpha}\right) |D|
<\infty,
$$
and hence is  compact.
 For $a\geq  0$, let
$ \{\lambda^{a, D}_k
: k=1, 2 \dots \}$ be the eigenvalues of $-(\Delta +a^\alpha
\Delta^{\alpha/2})|_D$, arranged in increasing order and repeated
according to multiplicity, and $\{\phi^{a, D}_k :k=1, 2, \dots\} $
be the corresponding eigenfunctions normalized to have unit
$L^2$-norm on $D$. Note that  $\{\phi^{a, D}_k :k=1, 2, \dots\} $
forms an orthonormal basis of $L^2(D; dx)$. It is well known that
when $a>0$, $\lambda^{a, D}_1$ is strictly positive and simple,
 and that $\phi^{a, D}_1  $ can be chosen to be
strictly positive on $D$. It follows from \cite[Theorem
1.1(ii)]{CS8} that the function $a\mapsto \lambda^{a, D}_1$ is
continuous
on $(0, M]$, and $\lim_{a\to 0+} \lambda^{a, D}_1 =
\lambda^{0, D}_1:= \min_{1\leq j\leq k} \lambda^{0, D_j}_1$,
where $D_1, \cdots, D_k$ are the connected components of $D$
and $\lambda^{0,D_j}_1$ is the first Dirichlet eigenvalue of
$-\Delta|_{D_j}$. Hence there is a constant $c_1=c_1(D, \alpha, M)\geq 1$
so that
\begin{equation}\label{e:4.1}
 1/c_1 \leq \lambda^{a, D}_1 \leq c_1
\qquad \hbox{for every } a\in (0, M].
\end{equation}
Using Sobolev embedding (see \cite[Example 5.1]{CS8}), it can be
shown that $\{\phi^{a, D}_1, \, a\in (0, M]\}$ is relatively compact
in $L^2(D; dx)$. Hence by  \cite[Theorem 1.1(ii)]{CS8} and the fact
that each $\lambda^{a, D}_1$ is simple for $a>0$, $a\to \phi^{a,
D}_1$ is continuous in $L^2(D; dx)$ in $a\in (0, M]$. Furthermore,
as $a\to 0+$,  any weak limit $\phi$ of $\phi^{a, D}_1$ is a unit
non-negative eigenfunction of $-\Delta|_D$ with eigenvalue
$\lambda^{0, D}_1$. Note that such $\phi$ may not be strictly
positive everywhere on $D$ when $D$ is disconnected; it is strictly
positive on least one component $D_j$ of $D$ where $\lambda^{0,
D_j}_1=\lambda^{0, D}_1$.
It follows that there is a constant
$c_2=c_2(D, \alpha, M)>0 $ so that
\begin{equation}\label{e:4.2}
\sup_{a\in (0, M]} \int_D \delta_D(x) \phi^{a, D}_1 (x) dx \leq c_2.
\end{equation}

Recall that $p^a_D (t, x, y)$ admits the following eigenfunction
expansion
$$
p^a_D (t, x, y) =\sum_{k=1}^\infty e^{-t \lambda^{a, D}_k} \phi^{a,
D}_k (x) \phi^{a, D}_k(y)
\qquad \hbox{for } t>0 \hbox{ and } x, y \in D.
$$
This implies that
\begin{equation}\label{e:4.3}
\int_{D\times D} \delta_D(x) p^a_D(t, x, y)  \delta_D(y) \, dx dy =
\sum_{k=1}^\infty e^{-t \lambda^{a, D}_k} \left( \int_D \delta_D(x)
\phi^{a, D}_k (x) dx\right)^2.
\end{equation}
Consequently, we have
\begin{equation}\label{e:4.4}
\int_{D\times D} \delta_D(x) p^a_D(t, x, y)  \delta_D(y) \, dx dy
\leq  e^{-t \lambda^{a, D}_1} \, \int_D \delta_D(x)^2 dx
\end{equation}
for all $a>0$ and $t>0$.
On the other hand, since
$$\phi^{a, D}_1 (x) = e^{\lambda^{a, D}_1} \int_D p^a_D (1, x, y)
\phi^{a, D}_1 (y) dy,
$$
by  the upper bound estimate in Theorem \ref{t:main}(i) and
\eqref{e:4.1}-\eqref{e:4.2} that there is a constant $c_3=c_3(D,
\alpha, M)>0 $ so that for every $a\in (0, M]$ and $x\in D$,
$$ \phi^{a, D}_1 (x)\leq c_3 \delta_D(x) \int_D \delta_D(y)
\phi^{a, D}_1 (y) dy \leq  c_2 c_3 \, \delta_D(x) .
$$
It now follows from \eqref{e:4.3}  that for every that for every
$a\in (0, M]$ and $t>0$
\begin{eqnarray}\label{e:4.5}
&&  \int_{D\times D}  \delta_D(x) p^a_D(t, x, y)  \delta_D(y) \, dx
dy
\geq  e^{-t \lambda^{a, D}_1} \, \left(\int_D \delta_D(x) \phi^{a,D}_1(x) dx
\right)^2 \nonumber  \\
&&\geq  e^{-t \lambda^{a, D}_1} \, \left(\int_D  (c_2c_3)^{-1}
\phi^{a,D}_1(x)^2 dx\right)^2 = (c_2c_3)^{-2} \, e^{-t \lambda^{a,
D}_1} .
\end{eqnarray}

\medskip

It suffices to prove (ii)-(iii) of Theorem \ref{t:main} for $T\geq
3$. For $t\geq T$ and $x, y\in D$, observe that
\begin{equation}\label{e:4.6}
p^a_D (t, x, y)=\int_{D\times D}
p^a_D (1, x, z) p^a_D(t-2, z, w)
 p^a_D(1, w, y) dz dw .
\end{equation}
Since $D$ is bounded,  we have by the upper bound estimate in
Theorem \ref{t:main}(i) and \eqref{e:4.4} that there are constants
$c_i=c_i(D, \alpha, M)>0$, $i=4, 5$ so that for every $a\in (0, M]$,
$t\geq T$ and $x, y\in D$,
\begin{equation}\label{e:4.7}
 p^a_D(t, x, y) \leq c_4 \delta_D (x)\delta_D(y)  \int_{D\times D}
 \delta_D(z) p^a_D (t-2, z, w) \delta_D(w) dz dw
 \leq   c_5 \, \delta_D(x) \delta_D (y) e^{-t\lambda^{a, D}_1}.
\end{equation}

(ii) Assume first that $D$ is connected. Since $D$ is bounded and
connected, we have by the lower bound estimate in Theorem
\ref{t:main}(i) and \eqref{e:4.5} that there are constants
$c_i=c_i(D, \alpha, M)>0$, $i=6, 7$,
so that for every $a\in (0, M]$, $t\geq T$ and $x, y\in D$,
\begin{equation}\label{e:4.8}
p^a_D(t, x, y) \geq  c_6\, \delta_D (x)\delta_D(y)
 \int_{D\times D} \delta_D(z) p^a_D (t-2, z, w) \delta_D(w) dz dw
 \geq c_7 \, \delta_D(x) \delta_D (y) e^{-t\lambda^{a, D}_1}.
\end{equation}
This combined with \eqref{e:4.7} proves Theorem \ref{t:main}(ii).
\medskip

(iii)  Now let consider the case that $D$ is disconnected.
Note that it follows from (ii) that for every $t\geq 1$, $x\in D$
and $y\in D(x)$,
\begin{equation}\label{e:4.8-2} p^a_D(t,x, y) \geq p^a_{D(x)}(t, x, y)
\geq c_{8} e^{-t \lambda_1^{a, D(x)}} \delta_D(x) \delta_D(y).
\end{equation}
Moreover, the above inequality, \eqref{e:4.7} and the two-sided
estimate in Theorem \ref{t:main}(i) yield that there are a constant
$c_{9}:=c_{9} (D, \alpha, M)\geq 1$  such that for every $a\in (0,
M]$, $t>0$ and $x\in D$,
\begin{equation}\label{e:4.9}
  c_{9}^{-1} \, e^{-t \lambda_1^{a, D(x)}}\delta_D(x)
  \leq \P_x  \left(
 \tau_{D(x)}^a >t\right) \leq c_{9} \, e^{-t \lambda_1^{a, D(x)}}\delta_D(x)
\end{equation}
and
\begin{equation}\label{e:4.10}
\P_x  (\tau_D^a >t) \leq  c_{9} \, e^{-t \lambda_1^{a,
D}}\delta_D(x).
\end{equation}
For $t\geq T$, $x\in D$ and $y\in D\setminus D(x)$, we have
by the boundedness of $D$,
\eqref{e:4.9} and the lower bound estimate in Theorem
\ref{t:main}(i) that
\begin{eqnarray}
p^a_D(t, x, y) &=& \E_x \left[ p^a_D\big(t-\tau^a_{D(x)},
X^a_{\tau^a_{D(x)}}, y\big); \,
\tau^a_{D(x)}<t \right] \nonumber \\
&=& \int_0^t \left( \int_{D(x)} p^a_{D(x)}(s, x, z)\left( \int_{D\setminus D(x)} J^a(z, w) p^a_D(t-s, w, y) dw \right)dz \right) ds \label{e:4.111} \\
&\ge & c_{10} \, a^\alpha \int_0^t \left( \int_{D(x)} p^a_{D(x)}(s,
x, z) dz \right) \left(
\int_{D(y)} p^a_{D(y)}(t-s, w, y) dw \right) ds \nonumber \\
&= & c_{10} \, a^\alpha \int_0^t \P_x \big(\tau^a_{D(x)}> s\big)
\, \P_y \big( \tau^a_{D(y)}>t-s \big) \, ds  \nonumber \\
&\ge& c_{10}c^{-2}_{9}\, a^\alpha \int_0^t e^{-s\, \lambda_1^{a,
D(x)}} \,
\delta_{D}(x) \, e^{-(t-s)  \lambda_1^{a, D(y)}}\,
\delta_D(y) \, ds
\label{e:4new1} \\
&\geq &
c_{10}c^{-2}_{9}\, a^\alpha t \,  e^{-t (\lambda_1^{a, D(x)}
\vee \lambda_1^{a, D(y)})} \, \delta_D(x) \delta_D(y).
\label{e:4new2}
\end{eqnarray}
On the other hand,
using \eqref{e:4.9}-\eqref{e:4.10},
we have from \eqref{e:4.111}
that for  $t\geq T$, $x\in D$ and $y\in D\setminus D(x)$,
\begin{eqnarray}
p^a_D(t, x, y) &\leq & c_{11} \, a^\alpha \int_0^t \left(
\int_{D(x)} p^a_{D(x)}(s, x, z) dz \right) \left(
\int_{D\setminus D(x)} p^a_D(t-s, w, y) dw \right) ds \nonumber \\
&\leq & c_{11} \, a^\alpha \int_0^t \P_x (\tau^a_{D(x)}> s)
\, \P_y ( \tau^a_D>t-s) ds  \nonumber \\
&\leq & c_{11}c^{2}_{9}\, a^\alpha \int_0^t e^{-s\, \lambda_1^{a,
D(x)}} \, \delta_{D}(x) \, e^{-(t-s)  \lambda_1^{a, D}}\,
\delta_D(y) ds
\nonumber \\
&= & c_{11}c^{2}_{9}\, a^\alpha \delta_D(x) \delta_D(y) e^{-t \,
\lambda_1^{a, D} }
\int_0^t e^{-s \left(\lambda_1^{a, D(x)}-\lambda_1^{a, D}\right) } ds \label{e:4.15} \\
&\leq &
c_{11}c^{2}_{9}\,  a^\alpha t\,  e^{-t \lambda_1^{a, D }} \,
\delta_D(x) \delta_D(y). \label{e:4.16}
\end{eqnarray}

Finally, by \eqref{e:4.7} and the same argument
that leads to
\eqref{e:4.15}, we have that for  $t\geq T$, $x, y\in D(x)$,
\begin{eqnarray}
&&p^a_D(t, x, y) \nonumber \\
&=& p^a_{D(x)}(t, x, y) + \E_x \left[ p^a_D\big(t-\tau^a_{D(x)},
X^a_{\tau^a_{D(x)}}, y\big); \,
\tau^a_{D(x)}<t \right] \nonumber \\
&\leq &  c_{12} \delta_D(x) \delta_D (y) e^{-t\lambda^{a, D(x)}_1} +
c_{12} \, a^\alpha \delta_D(x) \delta_D(y) e^{-t \, \lambda_1^{a, D}}
\int_0^t e^{-s \left(\lambda_1^{a, D(x)}-\lambda_1^{a, D}\right) } ds
\label{e:4.17}\\
&\le &   c_{12} \delta_D(x) \delta_D (y) e^{-t\lambda^{a, D(x)}_1} +
c_{12} \, a^\alpha\, t \, \delta_D(x)\delta_D(y) e^{-t \, \lambda_1^{a,
D} }.
\label{e:4.18}
\end{eqnarray}
Combining this with
\eqref{e:4.7}--\eqref{e:4.8-2},
\eqref{e:4new2} and \eqref{e:4.16}
completes the proof of Theorem \ref{t:main}(iii).
\qed

\medskip

\begin{remark}\label{R:4.1} \rm
In general, when passing from \eqref{e:4new1} to \eqref{e:4new2},
 from \eqref{e:4.15} to \eqref{e:4.16} and  from \eqref{e:4.17}
  to \eqref{e:4.18},
 a factor $t$ is needed in order to have the lower estimate and the upper
  estimate
  that is uniform in $a\in (0, M]$.
Note that for $D$ having at least two connected components,
$\lambda_1^{a, D(x)}>\lambda_1^{a, D}$ for every $a>0$. Since $D$ is
a bounded $C^{1,1}$ open set, it has only finite many connected
components $D_1, \cdots, D_k$. According to \cite[Theorem 1.1]{CS8},
as $a\to 0$, $\lambda_1^{a, D}$ converges to $\min_{1\leq j\leq k}
\lambda_1^{0, D_j}$, where $ \lambda_1^{0, D_j}$ is the first
Dirichlet eigenvalue of $-\Delta|_{D_j}$ on domain $D_j$. Let $j_0$
be such that $\lambda_1^{0, D_{j_0}}= \min_{1\leq j\leq k}
\lambda_1^{0, D_j}$. Then for every $x\in D_{j_0}$, we have
$\displaystyle  \inf_{a\in (0, M]} \left( \lambda_1^{a, D(x)}-
\lambda_1^{a, D}\right) =\lim_{a\to 0+} \lambda_1^{a, D(x)}-
\lambda_1^{a, D}=0$. Moreover, if $D$ has two connected components
$D_1$ and $D_2$ that are isometric to each other, then by
\cite[Theorem 1.1]{CS8}, for $x\in D_1$ and $y\in D_2$,
$$ \lim_{a\to 0+} \lambda^{a, D(x)}_1=\lambda^{0, D_1}_1
=\lambda^{0, D_2}_1=\lim_{a\to 0+} \lambda^{a, D(y)}_1.
$$
\qed
\end{remark}

\vskip 0.3truein

{\bf Zhen-Qing Chen}

Department of Mathematics, University of Washington, Seattle,
WA 98195, USA

E-mail: \texttt{zchen@math.washington.edu}

\bigskip

{\bf Panki Kim}

Department of Mathematical Sciences and Research Institute of Mathematics,
Seoul National University,
San56-1 Shinrim-dong Kwanak-gu,
Seoul 151-747, Republic of Korea

E-mail: \texttt{pkim@snu.ac.kr}

\bigskip

{\bf Renming Song}

Department of Mathematics, University of Illinois, Urbana, IL 61801, USA

E-mail: \texttt{rsong@math.uiuc.edu}

\begin{thebibliography}{99}

\bibitem{BBKRSV} K. Bogdan, T. Byczkowski, T. Kulczycki, M. Ryznar, R. Song and Z. Vondra\v{c}ek,
{\it Potential analysis of stable processes and its extesions.}
Lecture Notes in Math, {\bf 1980}, Springer, 2009.

\bibitem{BG} K. Bogdan and T. Grzywny,
Heat kernel of fractional Laplacian in cones. {\it Colloq. Math.} (to appear)


\bibitem{BGR} K. Bogdan, T. Grzywny and M. Ryznar,
Heat kernel estimates for the fractional Laplacian with Dirichlet
conditions. {\it Ann. Probab.} (to appear) arXiv:0905.2626v2.



\bibitem{CSS}
L. Caffarelli, S. Salsa and L. Silvestre,
Regularity estimates for the solution and the free boundary to the obstacle
problem for the fractional Laplacian.
{\it Invent. Math.} {\bf 171(1)} (2008) 425--461.

\bibitem{CaS} L. Caffarelli and L. Silvestre, Regularity theory for
fully nonlinear integro-differential equations.
{\it Comm. Pure Appl. Math. \bf 62} (2009),  597--638.


\bibitem{CV} L. Caffarelli and A. Vasseur,
Drift diffusion equations with fractional diffusion
and the quasi-geostrophic equation,  \ {\it  Ann.  Math.} (to appear).



\bibitem{C0} Z.-Q. Chen,
Multidimensional symmetric stable processes.  {\it Korean J. Comput.
Appl. Math.} {\bf 6} (1999), 227--266.


\bibitem{C} Z.-Q. Chen, Symmetric jump processes and their
heat kernel estimates. {\it  Sci. China Ser. A}, {\bf 52} (2009),
1423--1445.

\bibitem{CKS}  Z.-Q. Chen, P. Kim, and R. Song,
Heat kernel estimates for Dirichlet fractional Laplacian.
{\em J. European Math. Soc.}, (to appear).

\bibitem{CKS1} Z.-Q. Chen, P. Kim and R. Song,
Two-sided heat kernel estimates for censored stable-like processes.
{\em Probab. Theory Relat. Fields},
{\bf 146} (2010), 361--399.


\bibitem{CKS2} Z.-Q. Chen, P. Kim and R. Song, Sharp heat kernel
estimates for relativistic stable processes in open sets. Preprint,
2009,   arXiv:0908.1509

\bibitem{CKS3} Z.-Q. Chen, P. Kim and R. Song, Dirichlet Heat Kernel Estimates for
$\Delta^{\alpha /2}+ \Delta^{\beta /2}$. Preprint, 2009,
arXiv:0910.3266


\bibitem{CKSV1} Z.-Q. Chen, P. Kim, R. Song and Z. Vondracek,
Boundary Harnack principle for $\Delta + \Delta^{\alpha/2}$.
Preprint, 2009, arXiv:0908.1559v2

\bibitem{CKSV2} Z.-Q. Chen, P. Kim, R. Song and Z. Vondracek,
Sharp Green function estimates for $\Delta + \Delta^{\alpha/2}$ in
$C^{1, 1}$ open sets and their applications.
Preprint, 2009, arXiv:0912.2765

\bibitem{CK} Z.-Q. Chen and T. Kumagai,
Heat kernel estimates for stable-like processes on $d$-sets. {\em
Stoch. Proc. Appl.}   {\bf 108} (2003), 27--62.

\bibitem{CK2}  Z.-Q. Chen and  T. Kumagai,
Heat kernel estimates for jump processes of mixed types on metric
measure spaces. {\it Probab. Theory Relat. Fields},  {\bf 140}
(2008), 277--317.

\bibitem{CK08}  Z.-Q. Chen and  T. Kumagai,
A priori H\"older estimate, parabolic Harnack principle and heat
kernel estimates for diffusions with jumps. \  {\it Rev. Mat. Iberoam.}
to appear, 2009.

\bibitem{CS8} Z.-Q. Chen and R. Song, Continuity of eigenvalues for
subordinate processes in domains. {\it Math. Z.}, {\bf 252} (2006),
71--89.

\bibitem{CT}
Z.-Q.~Chen and J. Tokle,  Global heat kernel estimates for
fractional Laplacians in unbounded  open sets. {\em Probab. Theory
Relat. Fields}, DOI 10.1007/s00440-009-0256-0 (online first).

\bibitem{Ch}
S. Cho,
Two-sided global estimates of the Green's function of parabolic equations.
{\em Potential Analysis},
{\bf 25(4)}  (2006), 387--398.

\bibitem{D1}
E. B. Davies, Explicit constants for Gaussian upper bounds on heat
kernels. {\em Amer. J. Math.} {\bf 109} (1987), 319--333.

\bibitem{D2}
E. B. Davies, The equivalence of certain heat kernel and Green
function bounds. {\em J. Funct. Anal.} {\bf 71} (1987), 88--103.


\bibitem{D3}  E. B. Davies,
{\em Heat Kernels and Spectral Theory}. Cambridge University Press,
Cambridge, 1989.

\bibitem{DS} E. B. Davies and B. Simon, Ultracontractivity and heat kernels
for Schr\"odinger operator and Dirichlet Laplacians. {\it J. Funct.
Anal.}, {\bf 59} (1984), 335-395.

\bibitem{FOT} M. Fukushima, Y. Oshima and M. Takeda,
{\em Dirichlet Forms and Symmetric Markov Processes}. Walter De
Gruyter, Berlin, 1994.



\bibitem{KSZ}  J. Klafter, M.~F. Shlesinger and G. Zumofen,
Beyond Brownian motion.
{\it Physics Today, \bf 49} (1996), 33--39.


\bibitem{JW}  A.  Janicki and  A. Weron,
{\it Simulation and Chaotic Behavior of $\alpha$-Stable Processes}.
Dekker, 1994.




\bibitem{OS}
B. {\O}ksendal and A. Sulem, {\it Applied stochastic control of jump
diffusions},  2nd edition. Springer, Berlin, 2007.



\bibitem{SV07}
R.~Song and  Z.~Vondra\v{c}ek, Parabolic Harnack inequality for the
mixture of Brownian motion and stable process.
\ {Tohoku Math. J. (2)}, {\bf 59} (2007), 1--19.

\bibitem{SV08} R.~Song and Z.~Vondra\v{c}ek,  On the relationship between
subordinate killed and killed subordinate processes.
 \ {\it Elect. Commun. Probab.}  {\bf 13} (2008), 325--336.



\bibitem{Zq3} Q. S. Zhang, The boundary behavior of heat kernels of Dirichlet
Laplacians, {\em J. Differential Equations}, {\bf 182} (2002),
416--430.

\end{thebibliography}
\end{document}